\documentclass[10pt]{article}
\usepackage[utf8]{inputenc}
\usepackage{amsmath}
\usepackage{amssymb}
\usepackage{mathrsfs}
\usepackage{amsfonts}
\usepackage{mathrsfs,amscd,amssymb,amsthm,amsmath,bm,graphicx}
\usepackage{graphicx}
\usepackage{cite}
\usepackage{color}
\usepackage{enumitem}
\makeatletter

\renewcommand{\@seccntformat}[1]{{\csname the#1\endcsname}{\normalsize.}\hspace{.5em}}
\makeatother

\def \[{\begin{equation}}
\def \]{\end{equation}}

\newtheorem{thm}{Theorem}[section]

\begin{document}
\setlength{\baselineskip}{13pt}
\begin{center}{\Large \bf The expected values, variances and limiting distributions of Gutman index, Schultz index, multiplicative degree-Kirchhoff index and additive degree-Kirchhoff index for a class of random chain networks}

\vspace{4mm}

{\large Jia-Bao Liu $^{1,*}$, Qing Xie $^{1,*}$, Jiao-Jiao Gu $^{1}$}
\footnotetext{E-mail address: liujiabaoad@163.com, xieqingai25@163.com, gujiaojiaoajd@163.com.\\}
 \footnotetext{* Corresponding author.}
 \vspace{2mm}

{\small $^1$ School of Mathematics and Physics, Anhui Jianzhu University, Hefei 230601, P.R. China\\}
\vspace{2mm}
\end{center}

 {\noindent{\bf Abstract.}\ \
 There has been an upsurge of research on complex networks in recent years. The purpose of this paper is to study the mathematical properties of the random chain networks $PG_{n}$ with the help of graph theory. We first solve the expected value expressions of the Gutman index, Schultz index, multiplicative degree-Kirchhoff index and additive degree-Kirchhoff index, and then we get the explicit expression formulas of their variances. Finally, we find that their limiting distributions all have the probabilistic and statistical significance of normal distribution.

 \vspace{2mm}

\section{Introduction} \ \ \ \ \
Complex networks have become increasingly important in scientific research, especially in statistical physics and computer and information science \cite{1}, so researchers began to try to use this new theoretical tool to study various large-scale complex systems in real life, among them, the structure of complex system and the relationship between system structure and system function are hot topics. The application of graph theory, as the name implies, is to solve the corresponding problems in real life by combining the knowledge of graph theory, widely used in engineering design. Graph theory is a basic but very effective tool to study complex networks. There are a lot of complex systems in nature which can be described by network graphs. It is composed of many nodes and the edges between the nodes, where nodes represent different individuals, the edge represents the relationship between individuals. Topological index is one of the tools in graph theory, which can reflect the related properties and connected information of some kinds of topological structures. It has been extensively researched in graph theory. Therefore, this paper combines the three methods to measure the characteristics of complex networks quantitatively.

Suppose $G$ represents a simple undirected graph with edge set $E_{G}$ and vertex set $V_{G}$. In the graph $G$, the number of edges associated with a vertex $v$ is called the degree of the vertex, it is expressed by $d(v)$. The distance $d(u,v)$ between any two vertices $u$ and $v$ in the graph represents the length of the shortest path between them. To learn more about the relevant definitions, one can be referred to \cite{2}. The Wiener index $W(G)$ is the sum of the distances between all pairs of vertices in a graph $G$, denoted by
\begin{eqnarray*}
W(G)=\sum_{\{u,v\}\subseteq V_{G}}d(u,v).
\end{eqnarray*}

This invariant related to the distance of a graph was introduced into chemistry in 1947 \cite{3} and mathematics 30 years later \cite{4}. In 2016, Chen et al. \cite{5,6} studied the relationships between Wiener index and other distance-related topological indices in the tree-like polyphenyl systems. Shiu et al. \cite{7,8,9,10,11} calculated the Wiener indices about pericondensed hexagonal systems. At present, Wiener index is a kind of topological indices which is broadly used. For more applications of Wiener index in chemistry and mathematics, please refer to \cite{12,13,14,15}.

Based on the concept of distance and electrical network theory, the resistance distance $r(u,v)$ is derived, which was put forward by Klein and Randi\'{c}\cite{16}. When each edge in the graph is equivalent to a unit resistance, it can be considered as the effective resistance between two vertices. We can learn more details through \cite{17}.The Kirchhoff index $Kf(G)$, also known as full effective resistance, or effective graph resistance, is defined as the sum of the resistance distances of all vertex pairs in the graph \cite{18}, that is
\begin{eqnarray*}
Kf(G)=\sum_{\{u,v\} \subseteq V_{G}}r(u, v).
\end{eqnarray*}

Through the efforts of researchers, it has been found that the Kirchhoff index has played an important role in chemistry. For example, for fullerenes, linear hexagonal chains, and some special molecular graphs, such as distance regular graphs and M\"{o}bius ladders, we can evaluate the cyclicity of these polycyclic structures \cite{19,20,21,22,23}. At the same time, Kirchhoff index is also of great significance in mathematics. For the details of the mathematical properties, please refer to \cite{24,25,26,27,28,29,30}.

As for the Wiener index, the researchers improved it by weighting it. That is a graph $G=(V_{G}, E_{G})$ related to the weight function $w$: $V_{G}\longrightarrow {\mathbb{N}}^{+}$, expressed by $(G,w)$. Suppose $\oplus$ represents a calculation symbol in ${+}$,~${-}$,~$\times$,~$\div$, then the transformed Wiener index $W(G,w)$ can be expressed as
\begin{eqnarray}
W(G,w)=\frac{1}{2}\sum_{u \in V_{G}}\sum_{v \in V_{G}}\big(w(u) \oplus w(v)\big)d(u,v).
\end{eqnarray}

When $\oplus$ is the arithmetic operation $\times$ and set $w(\cdot) \equiv d(\cdot)$, the Gutman index is obtained. Then Eq. (1.1)  is equivalent to
\begin{eqnarray}
Gut(G)=\frac{1}{2}\sum_{u \in V_{G}}\sum_{v \in V_{G}}\big(d(u)d(v)\big)d(u,v)=\sum_{\{u,v\}\subseteq V_{G}}\big(d(u)d(v)\big)d(u,v).
\end{eqnarray}

When $\oplus$ is the arithmetic operation $+$ and set $w(\cdot) \equiv d(\cdot)$, the Schultz index is obtained. Then Eq. (1.1) is equivalent to
\begin{eqnarray}
S(G)=\frac{1}{2}\sum_{u \in V_{G}}\sum_{v \in V_{G}}\big(d(u)+d(v)\big)d(u,v)=\sum_{\{u,v\}\subseteq V_{G}}\big(d(u)+d(v)\big)d(u,v).
\end{eqnarray}

In the field of chemistry, Gutman index and Schultz index can accurately reflect the similar molecular structure characteristics as Wiener index. For more details, properties and applications of Gutman index and Schultz index, if you are interested, you can refer to the article\cite{31,32,33,34}.

Similarly, the researchers also improved and optimized the Kirchhoff index. Considering the degree of vertices in the graph, the multiplicative degree-Kirchhoff index was put forward in 2007 by Chen and Zhang\cite{35}, which was expressed as
\begin{eqnarray}
Kf^{*}(G)=\sum_{\{u,v\} \subseteq V_{G}}d(u)d(v)r(u, v).
\end{eqnarray}
In 2012, the additive degree-Kirchhoff index was proposed by Gutman, Feng and Yu \cite{36}, denoted by
\begin{eqnarray}
Kf^{+}(G)=\sum_{\{u,v\} \subseteq V_{G}}\big(d(u)+d(v)\big)r(u, v).
\end{eqnarray}

Kirchhoff index and degree Kirchhoff index are used as graph parameters based on resistance distance in physics and chemistry. And network science also has a wide range of applications. For an electrical network, the smaller the Kirchhoff index is, the less power the network consumes per unit time. Chemically, Kirchhoff index can be used to describe the structural characteristics of molecules and define the topological radius of polymers. In network science, the smaller the Kirchhoff index is, the stronger the robustness of the network is.  These fields have attracted the attention of more and more researchers. Readers can refer to the reference \cite{37,38} for more information.

The random chain networks graphs studied in this paper are finite 2-connected graphs, which can be regarded as composed of $n$ pentagons, which are gradually increased by connecting edges at the ends of the pentagons. That is, $PG_{n}$ is obtained by $PG_{n-1}$ and then randomly connecting a pentagon at the end, as shown in figure 2. When $n$ is equal to 1, 2, 3, as shown in figure 1. For $n \geq 3$, there are two ways to add pentagons at the end, and the results can be expressed as $PG_{n}^{1}$ and $PG_{n}^{2}$, as shown in figure 3. For such a random chain networks, any step for $k=3, 4, \cdot \cdot \cdot, n$ is stochastic, and their probabilities are $p_{1}$ and $p_{2}$, respectively:

$\bullet$ (i) $PG_{n-1}\longrightarrow PG^{1}_{n}$ with probability $p_{1}$,

$\bullet$ (ii) $PG_{n-1}\longrightarrow PG^{2}_{n}$ with probability $p_{2}=1-p_{1}$,
\\\text{in which the probabilities $p_{1}$ and $p_{2}$ are constants and independent to the step $k$ at the same time.}\par

Motivated by \cite{39}, we use two random variables $Z^{1}_{n}$ and $Z^{2}_{n}$ to represent our choices. If our choice is $PG^{i}_{n}$, we put $Z^{i}_{n}=1$, otherwise $Z^{i}_{n}=0$ ($i=1,~ 2$). One holds that
\begin{eqnarray}
\mathbb{P}(Z^{i}_{n}=1)=p_{i},~~~\mathbb{P}(Z^{i}_{n}=0)=1-p_{i},~~~i=1,~2,
\end{eqnarray}
and $Z^{1}_{n}+Z^{2}_{n}=1$.

\begin{figure}[htbp]
\centering\includegraphics[width=10.92cm,height=6.136cm]{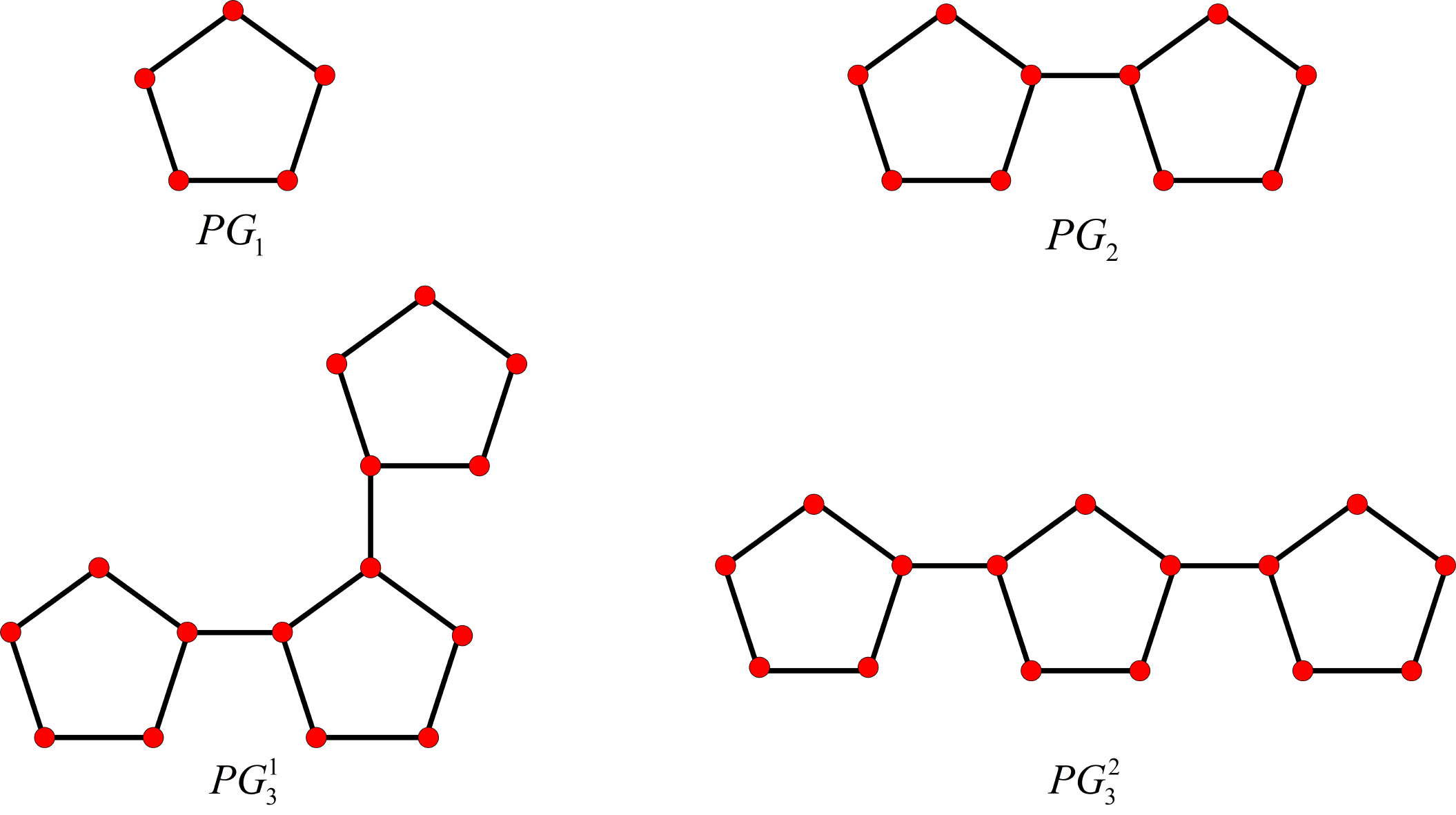}
\caption{ Three kinds of random chain networks.}
\end{figure}

\begin{figure}[htbp]
\centering\includegraphics[width=5.88cm,height=3.304cm]{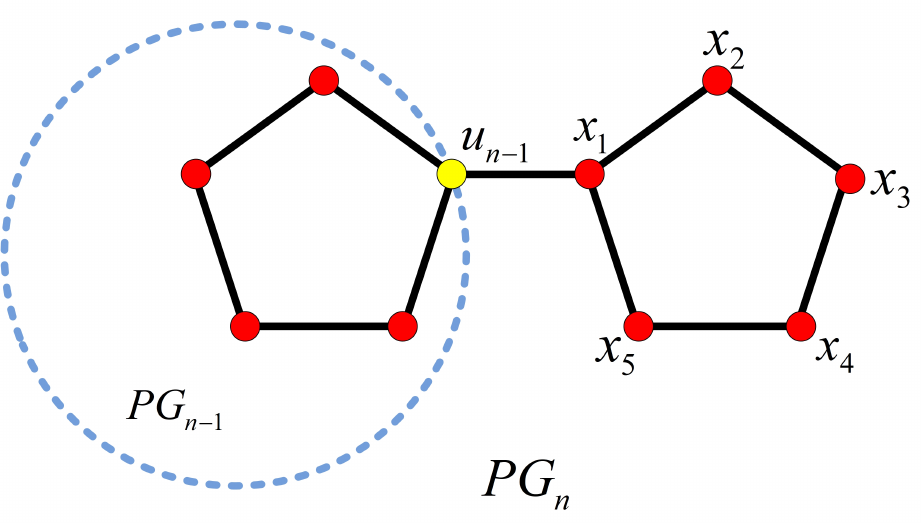}
\caption{ A random chain networks $PG_{n}$.}
\end{figure}

Yang and Zhang\cite{40} and Ma et al.\cite{41} determined explicit formulas of $W(G)$ and $\mathbb{E}\big(W(G)\big)$ for a class of random chain networks composed of hexagons, respectively. Similarly, Huang, Kuang and Deng \cite{42} obtained $\mathbb{E}\big(Kf(G)\big)$ for a class of helical random chain networks composed of hexagons. Wei and Shiu \cite{43} put forward the expression of $\mathbb{E}\big(W(G)\big)$ for random polygon chain networks and proved the asymptotic property of its expected values. At the same way, Zhang, Li, Li and Zhang \cite{44} obtained the simple formulas of $\mathbb{E}\big(Gut(G_{n})\big)$, $\mathbb{E}\big(S(G_{n})\big)$, $\mathbb{E}\big(Kf^{*}(G_{n})\big)$ and $\mathbb{E}\big(Kf^{+}(G_{n})\big)$ for a class of random chain networks.

\begin{figure}[htbp]
\centering\includegraphics[width=12.936cm,height=4.8048cm]{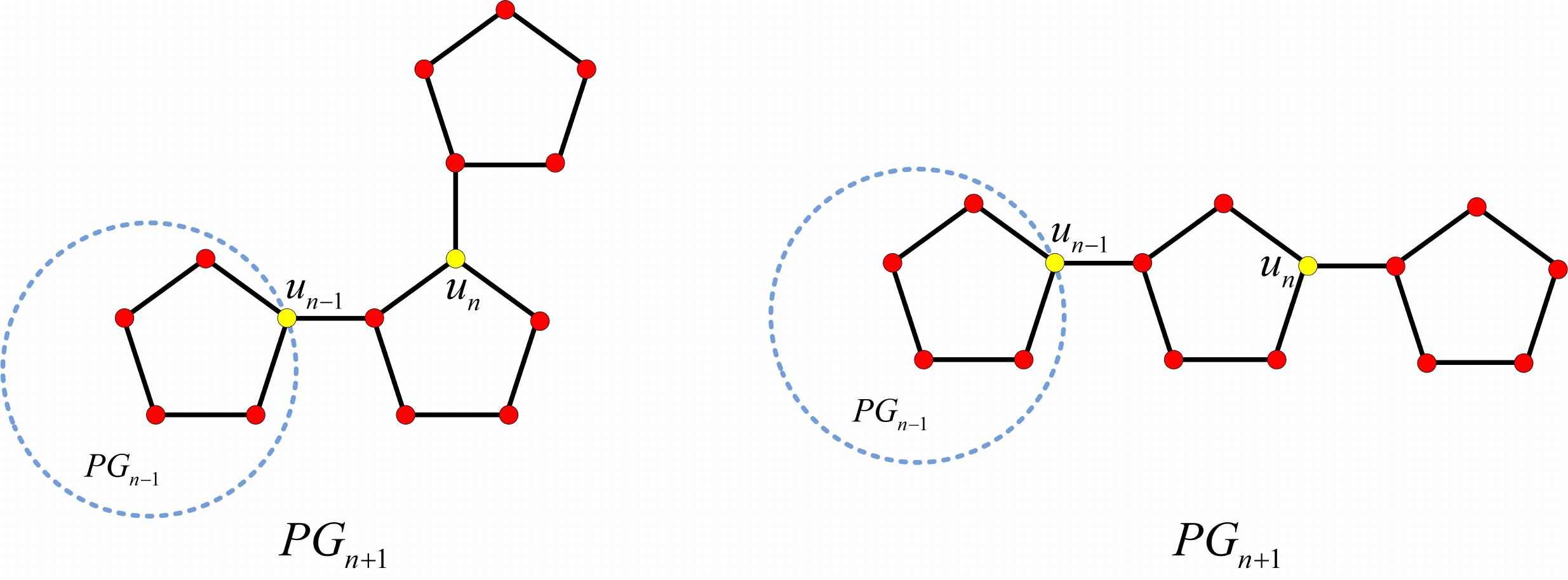}
\caption{ Two connection modes in chain networks.}
\end{figure}

Motivated by \cite{44,45,46,47,48}, for this article, in section 2, we determine the explicit formulas of $\mathbb{E}\big(Gut(PG_{n})\big)$, $\mathbb{E}\big(S(PG_{n})\big)$, $\mathbb{E}\big(Kf^{*}(PG_{n})\big)$ and $\mathbb{E}\big(Kf^{+}(PG_{n})\big)$ for the random chain networks. In section 3, we obtain explicit formulas of $Var\big(Gut(PG_{n})\big)$, $Var\big(S(PG_{n})\big)$, $Var\big(Kf^{*}(PG_{n})\big)$ and $Var\big(Kf^{+}(PG_{n})\big)$ based on the known results of others’ paper, and prove that these four indices of the random chain networks asymptotically obey normal distributions.

For the purposes of this paper, we must ensure that the following Hypotheses hold.

\newtheorem{hyp}{Hypothesis}
\begin{hyp}
It is randomly and independently to choice a way attaching the new terminal pentagon $O_{n+1}$ to $PG_{n}$, $n=2, 3, \cdot \cdot \cdot$. To be more precisely, the sequences of random variables $\{Z^{1}_{n}, Z^{2}_{n}\}^{\infty}_{n=2}$ are independently and must satisfy Eq. (1.6).
\end{hyp}
\begin{hyp}
For $i\in \{1, 2\}$, we put $0< p_{i} < 1$.
\end{hyp}
Under the conditions that Hypotheses 1 and 2.

(a) The analytical expressions of the variances of $Gut(PG_{n})$, $S(PG_{n})$, $Kf^{*}(PG_{n})$ and $Kf^{+}(PG_{n})$ are obtained;

(b) When $n \rightarrow \infty$, we verify that the random variables $Gut(PG_{n})$, $S(PG_{n})$, $Kf^{*}(PG_{n})$ and $Kf^{+}(PG_{n})$ asymptotically obey normal distributions. It is evident to see that
\begin{eqnarray*}
\lim_{n\rightarrow \infty}\sup_{a\in \mathbb{R}}\mid \mathbb{P}\big(\frac{X_{n}-\mathbb{E}(X_{n})}{\sqrt{Var(X_{n})}}\leq a\big)-{\int}^{a}_{-\infty}\frac{1}{\sqrt{2\pi}}e^{-\frac{t^{2}}{2}}dt \mid =0,
\end{eqnarray*}
where $\mathbb{E}(X_{n})$ and $Var(X_{n})$ represent the expectation and variance of this random variable $X_{n}$ respectively.

In this paper, assume that $f(x)$ and $h(x)$ are two functions of $x$. We put $f(x) \asymp h(x)$ if $\lim\limits_{x \rightarrow \infty}\frac{f(x)}{h(x)}=1$, and put $f(x)=O\big(h(x)\big)$ if $\lim\limits_{x \rightarrow \infty}\frac{f(x)}{h(x)}=0$.

\section{The expected values of $Gut(PG_{n})$, $S(PG_{n})$, $Kf^{*}(PG_{n})$ and $Kf^{+}(PG_{n})$ of the random chain networks}
\ \ \ \
For the random chain networks $PG_{n}$, we find that $Gut(PG_{n})$, $S(PG_{n})$, $Kf^{*}(PG_{n})$ and $Kf^{+}(PG_{n})$ are random variables. Then, we determine the analytic expressions of $\mathbb{E}\big(Gut(PG_{n})\big)$, $\mathbb{E}\big(S(PG_{n})\big)$, $\mathbb{E}\big(Kf^{*}(PG_{n})\big)$ and $\mathbb{E}\big(Kf^{+}(PG_{n})\big)$ in this section.

In fact, $PG_{n}$ is organized by adding a new terminal pentagon $O_{n}$ to $PG_{n-1}$ by an edge, where the vertices of $O_{n}$ are labelled as $x_{1}, ~x_{2},~ x_{3}, ~x_{4},~ x_{5}$ in clockwise direction. For all $v \in V_{PG_{n-1}}$, one has
\begin{equation}
\begin{split}
&d(x_{1},v)=d(u_{n-1},v)+1,~~~d(x_{2},v)=d(u_{n-1},v)+2,~~~d(x_{3},v)=d(u_{n-1},v)+3,\\
&d(x_{4},v)=d(u_{n-1},v)+3,~~~d(x_{5},v)=d(u_{n-1},v)+2,\\
&\sum_{v \in V_{PG_{n-1}}}d_{PG_{n}}(v)=12n-13, ~~and ~~\sum_{v \in V_{PG_{n}}}d_{PG_{n+1}}(v)=12n-1.
\end{split}
\end{equation}
Meanwhile,
\begin{equation}
\begin{split}
&\sum_{i=1}^{5}d(x_{i})d(x_{1},x_{i})=22,~~~\sum_{i=1}^{5}d(x_{i})d(x_{2},x_{i})=13,~~~\sum_{i=1}^{5}d(x_{i})d(x_{3},x_{i})=14,\\
&\sum_{i=1}^{5}d(x_{i})d(x_{4},x_{i})=15,~~~\sum_{i=1}^{5}d(x_{i})d(x_{5},x_{i})=13.\\
\end{split}
\end{equation}

\begin{thm}
For $n\geq 1$, the analytic expression of $\mathbb{E}\big(Gut(PG_{n})\big)$ for the random chain networks is
\begin{eqnarray*}
\mathbb{E}\big(Gut(PG_{n})\big)&=&(72-24p_{1})n^{3}+(72p_{1}-12)n^{2}+(1-48p_{1})n-1.
\end{eqnarray*}
\end{thm}

\noindent{\bf Proof of Theorem 2.1.}  By Eq. (1.2), one can be convinced that
\begin{eqnarray*}
Gut(PG_{n+1})&=&\sum_{\{u,v\}\subseteq V_{PG_{n}}}d(u)d(v)d(u,v)+\sum_{v\in V_{PG_{n}}}\sum_{x_{i}\in V_{O_{n+1}}}d(v)d(x_{i})d(v,x_{i})\\
&&+\sum_{\{x_{i},x_{j}\}\subseteq V_{O_{n+1}}}d(x_{i})d(x_{j})d(x_{i},x_{j}).
\end{eqnarray*}
Note that
\begin{align*}
&\sum_{\{u,v\}\subseteq V_{PG_{n}}}d(u)d(v)d(u,v)\\
&=\sum_{\{u,v\}\subseteq V_{PG_{n}}\setminus \{u_{n}\}}d(u)d(v)d(u,v)+\sum_{v\in V_{PG_{n}}\setminus \{u_{n}\}}d_{PG_{n+1}}(u_{n})d(v)d(u_{n},v)\\
&=\sum_{\{u,v\}\subseteq V_{PG_{n}}\setminus \{u_{n}\}}d(u)d(v)d(u,v)+\sum_{v\in V_{PG_{n}}\setminus \{u_{n}\}}\big(d_{PG_{n}}(u_{n})+1\big)d(v)d(u_{n},v)\\
&=Gut(PG_{n})+\sum_{v\in V_{PG_{n}}}d(v)d(u_{n},v).
\end{align*}
Recall that $d(x_{1})=3$ and $d(x_{i})=2$ for $i \in \{2, ~3,~ 4,~ 5\}$. From Eq. (2.7), we have
\begin{eqnarray*}
\sum_{v\in V_{PG_{n}}}\sum_{x_{i}\in V_{O_{n+1}}}d(v)d(x_{i})d(v,x_{i})&=&\sum_{v\in V_{PG_{n}}}d(v)\bigg[3\big(d(u_{n},v)+1\big)+2\big(d(u_{n},v)+2\big)+2\big(d(u_{n},v)+3\big)\\
&&+2\big(d(u_{n},v)+3\big)+2\big(d(u_{n},v)+2\big)\bigg]\\
&=&\sum_{v\in V_{PG_{n}}}d(v)\big(11d(u_{n},v)+23\big)\\
&=&11\sum_{v\in V_{PG_{n}}}d(v)d(u_{n},v)+276n-23.
\end{eqnarray*}
From Eq. (2.8), one follows that
\begin{eqnarray*}
\sum_{\{x_{i},x_{j}\}\subseteq V_{O_{n+1}}}d(x_{i})d(x_{j})d(x_{i},x_{j})&=&\frac{1}{2}\sum^{5}_{i=1}d(x_{i})\big(\sum^{5}_{j=1}d(x_{j})d(x_{i},x_{j})\big)\\
&=&\frac{1}{2}(3\times12+2\times 2\times 13+2\times 2\times 14)\\
&=&72.
\end{eqnarray*}
Then
\begin{eqnarray}
Gut(PG_{n+1})=Gut(PG_{n})+12\sum_{v\in V_{PG_{n}}}d(v)d(u_{n},v)+276n+49.
\end{eqnarray}

For the random chain networks $PG_{n}$, we obtain that $\sum_{v\in V_{PG_{n}}}d(v)d(u_{n},v)$ is a stochastic variate. Let
\begin{eqnarray*}
A_{n}:=\mathbb{E}\big(\sum_{v\in V_{PG_{n}}}d(v)d(u_{n},v)\big).
\end{eqnarray*}
By using above formula and Eq. (2.9), we can get the following relation for $\mathbb{E}\big(Gut(PG_{n+1})\big)$ of the random chain networks. One sees

\begin{eqnarray}
\mathbb{E}\big(Gut(PG_{n+1})\big)=\mathbb{E}\big(Gut(PG_{n})\big)+12A_{n}+276n+49.
\end{eqnarray}

Then, we go on to consider the following two possible cases.

{\bf Case 1.} $PG_{n}\longrightarrow PG^{1}_{n+1}$.

In this case, $u_{n}$ (of $PG_{n}$) overlaps with $x_{2}$ or $x_{5}$ (of $O_{n}$). Therefore, $\sum\limits_{v\in V_{PG_{n}}}d(v)d(u_{n},v)$ is rewritten as $\sum\limits_{v\in V_{PG_{n}}}d(v)d(x_{2},v)$ or $\sum\limits_{v\in V_{PG_{n}}}d(v)d(x_{5},v)$ with probability $p_{1}$.

{\bf Case 2.} $PG_{n}\longrightarrow PG^{2}_{n+1}$.

In this case, $u_{n}$ (of $PG_{n}$) overlaps with $x_{3}$ or $x_{4}$ (of $O_{n}$). Therefore, $\sum\limits_{v\in V_{PG_{n}}}d(v)d(u_{n},v)$ is rewritten as $\sum\limits_{v\in V_{PG_{n}}}d(v)d(x_{3},v)$ or $\sum\limits_{v\in V_{PG_{n}}}d(v)d(x_{4},v)$ with probability $p_{2}=1-p_{1}$.

Together with the above Cases, by applying the expectation operator and Eq. (2.7), one follows that
\begin{eqnarray*}
A_{n}&=&p_{1}\sum\limits_{v\in V_{PG_{n}}}d(v)d(x_{2},v)+(1-p_{1})\sum\limits_{v\in V_{PG_{n}}}d(v)d(x_{3},v),
\end{eqnarray*}
where
\begin{eqnarray*}
\sum\limits_{v\in V_{PG_{n}}}d(v)d(x_{2},v)
&=&\sum\limits_{v\in V_{PG_{n-1}}}d(v)\big(d(u_{n-1},v)+2\big)+\sum\limits_{v\in V_{O_{n}}}d(v)d(x_{2},v)\\
&=&\sum\limits_{v\in V_{PG_{n-1}}}d(v)d(u_{n-1},v)+2\sum\limits_{v\in V_{PG_{n-1}}}d(v)+13,\\
\end{eqnarray*}
\begin{eqnarray*}
\sum\limits_{v\in V_{PG_{n}}}d(v)d(x_{3},v)
&=&\sum\limits_{v\in
V_{PG_{n-1}}}d(v)\big(d(u_{n-1},v)+3\big)+\sum\limits_{v\in V_{O_{n}}}d(v)d(x_{3},v)\\
&=&\sum\limits_{v\in V_{PG_{n-1}}}d(v)d(u_{n-1},v)+3\sum\limits_{v\in V_{PG_{n-1}}}d(v)+14.
\end{eqnarray*}
Then, we obtain
\begin{eqnarray*}
A_{n}&=&p_{1}\big[\sum\limits_{v\in V_{PG_{n-1}}}d(v)d(u_{n-1},v)+2\sum\limits_{v\in V_{PG_{n-1}}}d(v)+13\big]\\
&&+(1-p_{1})\big[\sum\limits_{v\in V_{PG_{n-1}}}d(v)d(u_{n-1},v)+3\sum\limits_{v\in V_{PG_{n-1}}}d(v)+14\big]\\
&=&p_{1}\big[A_{n-1}+2(12n-13)+13\big]+(1-p_{1})\big[A_{n-1}+3(12n-13)+14\big]\\
&=&A_{n-1}+(36-12p_{1})n+(12p_{1}-25).
\end{eqnarray*}
Meanwhile, for $n=1$, the boundary condition is
\begin{eqnarray*}
A_{1}=\sum_{v\in V_{PG_{1}}}d(v)d(u_{1},v)=12.
\end{eqnarray*}
Using above condition and the recurrence relation with respect to $A_{n}$, it is no hard to obtain
\begin{eqnarray}
A_{n}&=&(18-6p_{1})n^{2}+(6p_{1}-7)n+1.
\end{eqnarray}
From Eq. (2.10), it holds that
\begin{eqnarray*}
\mathbb{E}\big(Gut(PG_{n+1})\big)&=&\mathbb{E}\big(Gut(PG_{n})\big)+12A_{n}+276n+49\\
&=&\mathbb{E}\big(Gut(PG_{n})\big)+12\big((18-6p_{1})n^{2}+(6p_{1}-7)n+1\big)+276n+49.
\end{eqnarray*}
For $n=1$, we obtain $\mathbb{E}\big(Gut(PG_{1})\big)=60$.
Similarly, according to the recurrence relation related to $\mathbb{E}\big(Gut(PG_{n})\big)$, we have
\begin{eqnarray*}
\mathbb{E}\big(Gut(PG_{n})\big)&=&(72-24p_{1})n^{3}+(72p_{1}-12)n^{2}+(1-48p_{1})n-1,
\end{eqnarray*}
as desired.\hfill\rule{1ex}{1ex}\

\begin{thm}
For $n\geq 1$, the analytic expression of $\mathbb{E}\big(S(PG_{n})\big)$ for the random chain networks is
\begin{eqnarray*}
\mathbb{E}\big(S(PG_{n})\big)&=&(60-20p_{1})n^{3}+(60p_{1}+7)n^{2}-(40p_{1}+7)n.
\end{eqnarray*}
\end{thm}

\noindent{\bf Proof of Theorem 2.2.} Notice that the random chain networks $PG_{n+1}$ is organized by adding a new terminal pentagon $O_{n+1}$ to $PG_{n}$ by an edge. By Eq. (1.3), one has
\begin{eqnarray}
S(PG_{n+1})=\Delta_{1}+\Delta_{2}+\Delta_{3},
\end{eqnarray}
where
\begin{eqnarray*}
\Delta_{1}=\sum_{\{u,v\}\subseteq V_{PG_{n}}}\big(d(u)+d(v)\big)d(u,v),\\
\Delta_{2}=\sum_{v\in V_{PG_{n}}}\sum_{x_{i}\in V_{O_{n+1}}}\big(d(v)+d(x_{i})\big)d(v,x_{i}),\\
\Delta_{3}=\sum_{\{x_{i},x_{j}\}\subseteq V_{O_{n+1}}}\big(d(x_{i})+d(x_{j})\big)d(x_{i},x_{j}).\\
\end{eqnarray*}
Let $d(u_{n}\mid PG_{n}):=\sum\limits_{v\in V_{PG_{n}}}d(u_{n},v)$. It is routine to check that
\begin{eqnarray*}
\Delta_{1}&=&\sum_{\{u,v\}\subseteq V_{PG_{n}}\setminus \{u_{n}\}}\big(d(u)+d(v)\big)d(u,v)+\sum_{v\in V_{PG_{n}}\setminus \{u_{n}\}}\big(d_{PG_{n+1}}(u_{n})+d(v)\big)d(u_{n},v)\\
&=&\sum_{\{u,v\}\subseteq V_{PG_{n}}\setminus \{u_{n}\}}\big(d(u)+d(v)\big)d(u,v)+\sum_{v\in V_{PG_{n}}\setminus \{u_{n}\}}\big(d_{PG_{n}}(u_{n})+1+d(v)\big)d(u_{n},v)\\
&=&S(PG_{n})+d(u_{n}\mid PG_{n}).
\end{eqnarray*}
Note that $PG_{n}$ has $5n$ vertices. Note that $d(x_{1})=3$, $d(x_{i})=2$ for $i \in \{2,~ 3, ~4, ~5\}$. By Eq. (2.7), we can know
\begin{eqnarray*}
\Delta_{2}&=&\sum_{v\in V_{PG_{n}}}\sum_{x_{i}\in V_{O_{n+1}}}d(v)d(v,x_{i})+\sum_{v\in V_{PG_{n}}}\sum_{x_{i}\in V_{O_{n+1}}}d(x_{i})d(x_{i},v)\\
&=&\sum_{v\in V_{PG_{n}}}d(v)\big(\sum_{x_{i}\in V_{O_{n+1}}}d(v,x_{i})\big)+\sum_{v\in V_{PG_{n}}}\big(\sum_{x_{i}\in V_{O_{n+1}}}d(x_{i})d(v,x_{i})\big)\\
&=&\sum_{v\in V_{PG_{n}}}d(v)\big(5d(u_{n},v)+11\big)+\sum_{v\in V_{PG_{n}}}\big(11d(u_{n}+23)\big)\\
&=&5\sum_{v\in V_{PG_{n}}}d(v)d(u_{n},v)+11(12n-1)+11\sum_{v\in V_{PG_{n}}}d(u_{n},v)+23\times 5n\\
&=&5\sum_{v\in V_{PG_{n}}}d(v)d(u_{n},v)+11dis(u_{n}\mid PG_{n})+247n-11.
\end{eqnarray*}
Note that $\sum\limits_{i=1}^{5}d(x_{k},x_{i})=6$ for $k=1, ~2,\cdot \cdot \cdot, ~5$. From Eq. (2.8), one sees that
\begin{eqnarray*}
\Delta_{3}&=&\sum_{\{x_{i},x_{j}\}\subseteq V_{O_{n+1}}}\big(d(x_{i})+d(x_{j})\big)d(x_{i},x_{j})=\frac{1}{2}\sum^{5}_{i=1}\sum^{5}_{j=1}\big(d(x_{i})+d(x_{j})\big)d(x_{i},x_{j})\\
&=&\sum^{5}_{i=1}\sum^{5}_{j=1}d(x_{i})d(x_{i},x_{j})=6\times(3+2\times 4)=66.
\end{eqnarray*}
Then, Eq. (2.12) can be rewritten as
\begin{eqnarray}
S(PG_{n+1})=S(PG_{n})+5\sum_{v\in V_{PG_{n}}}d(v)d(u_{n},v)+12d(u_{n}\mid PG_{n})+247n+55.
\end{eqnarray}

For the random chain networks $PG_{n}$, we know that $d(u_{n}\mid PG_{n})$ is a random variable. Let
\begin{eqnarray*}
B_{n}:=\mathbb{E}\big(d(u_{n}\mid PG_{n})\big).
\end{eqnarray*}
According to above formula and Eq. (2.13), we have the following relation for $\mathbb{E}\big(S(PG_{n+1})\big)$ of the random chain networks. It holds that
\begin{eqnarray}
\mathbb{E}\big(S(PG_{n+1})\big)=\mathbb{E}\big(S(PG_{n})\big)+5A_{n}+12B_{n}+247n+55.
\end{eqnarray}
We proceed by taking into account the following two cases.

{\bf Case 1.} $PG_{n}\longrightarrow PG^{1}_{n+1}$.

In this case, $u_{n}$ (of $PG_{n}$) coincides with $x_{2}$ or $x_{5}$ (of $O_{n}$). Then, $d(u_{n}\mid PG_{n})$ is given by $d(x_{2}\mid COC_{n})$ or $d(x_{5}\mid PG_{n})$ with probability $p_{1}$.

{\bf Case 2.} $PG_{n}\longrightarrow PG^{2}_{n+1}$.

In this case, $u_{n}$ (of $PG_{n}$) coincides with $x_{3}$ or $x_{4}$ (of $O_{n}$). Then, $d(u_{n}\mid PG_{n})$ is given by $d(x_{3}\mid PG_{n})$ or $d(x_{4}\mid PG_{n})$ with probability $p_{2}=1-p_{1}$.

Together with Cases 1 and 2, we obtain $B_{n}$ equals to
\begin{align*}
p_{1}\cdot d(x_{2}\mid PG_{n})+(1-p_{1})\cdot d(x_{3}\mid PG_{n}),
\end{align*}
where
\begin{eqnarray*}
&d(x_{2}\mid PG_{n})=d(u_{n-1}\mid PG_{n-1})+2\times5(n-1)+6=d(u_{n-1}\mid PG_{n-1})+10n-4,\\
&d(x_{3}\mid PG_{n})=d(u_{n-1}\mid PG_{n-1})+3\times5(n-1)+6=d(u_{n-1}\mid PG_{n-1})+15n-9.
\end{eqnarray*}
Then, as an immediate consequence, we have
\begin{eqnarray*}
B_{n}&=&p_{1}\big(B_{n-1}+10n-4\big)+(1-p_{1})\big(B_{n-1}+15n-9\big)\\
&=&B_{n-1}+(15-5p_{1})n+(5p_{1}-9).
\end{eqnarray*}
When $n=1$, we can easily find that
\begin{eqnarray*}
B_{1}=\mathbb{E}\big(d(u_{1}\mid PG_{1})\big)=6.
\end{eqnarray*}
Then, using above formula and the recurrence relation, it is routine to check that
\begin{eqnarray*}
B_{n}&=&(\frac{15}{2}-\frac{5}{2}p_{1})n^{2}+(\frac{5}{2}p_{1}-\frac{3}{2})n.
\end{eqnarray*}
From Eq. (2.11), we have
\begin{eqnarray*}
A_{n}&=&(18-6p_{1})n^{2}+(6p_{1}-7)n+1.
\end{eqnarray*}
By using Eq. (2.14), we get
\begin{eqnarray*}
\mathbb{E}\big(S(PG_{n+1})\big)&=&\mathbb{E}\big(S(PG_{n})\big)+5A_{n}+12B_{n}+247n+55\\
&=&\mathbb{E}\big(S(PG_{n})\big)+5\big((18-6p_{1})n^{2}+(6p_{1}-7)n+1\big)\\
&&+12\big((7.5-2.5p_{1})n^{2}+(2.5p_{1}-1.5)n\big)+247n+55\\
&=&\mathbb{E}\big(S(PG_{n})\big)+(180-60p_{1})n^{2}+(60p_{1}+194)n+60.
\end{eqnarray*}
As an immediate consequence, we find $\mathbb{E}\big(S(PG_{1})\big)=60$. Then, we arrive at
\begin{eqnarray*}
\mathbb{E}\big(S(PG_{n})\big)&=&(60-20p_{1})n^{3}+(60p_{1}+7)n^{2}-(40p_{1}+7)n,
\end{eqnarray*}
as desired.\hfill\rule{1ex}{1ex}\

Then, we determine the analytic expressions of $\mathbb{E}\big(Kf^{*}(PG_{n})\big)$ and $\mathbb{E}\big(Kf^{+}(PG_{n})\big)$ now. Recall, $PG_{n+1}$ is organized by adding a new terminal pentagon $O_{n+1}$ to $PG_{n}$ by an edge, where the vertices of $O_{n+1}$ are labelled as $x_{1}, ~x_{2},~ x_{3}, ~x_{4},~ x_{5}$ in clockwise direction. For all $v \in V_{PG_{n}}$, one has
\begin{equation}
\begin{split}
&r(x_{1},v)=r(u_{n},v)+1,~~~~~~~~~r(x_{2},v)=r(u_{n},v)+1+\frac{4}{5},~~~r(x_{3},v)=r(u_{n},v)+1+\frac{6}{5},\\
&r(x_{4},v)=r(u_{n},v)+1+\frac{6}{5},~~~r(x_{5},v)=r(u_{n},v)+1+\frac{4}{5},\\
&\sum_{v \in V_{PG_{n}}}d_{PG_{n+1}}(v)=12n-1.
\end{split}
\end{equation}
Therefore,
\begin{equation}
\begin{split}
&\sum_{i=1}^{5}d(x_{i})r(x_{1},x_{i})=8,~~~~~\sum_{i=1}^{5}d(x_{i})r(x_{2},x_{i})=\frac{44}{5},~~~\sum_{i=1}^{5}d(x_{i})r(x_{3},x_{i})=\frac{46}{5},\\
&\sum_{i=1}^{5}d(x_{i})r(x_{4},x_{i})=\frac{46}{5},~~~\sum_{i=1}^{5}d(x_{i})r(x_{5},x_{i})=\frac{44}{5}.\\
\end{split}
\end{equation}

\begin{thm}
For $n\geq 1$, the analytic expression of $\mathbb{E}\big(Kf^{*}(PG_{n})\big)$ for the random chain networks is
\begin{eqnarray*}
\mathbb{E}\big(Kf^{*}(PG_{n})\big)&=&(\frac{264}{5}-\frac{48}{5}p_{1})n^{3}+(\frac{144}{5}p_{1}-\frac{12}{5})n^{2}+(\frac{193}{5}-\frac{96}{5}p_{1})n-49.
\end{eqnarray*}
\end{thm}

\noindent{\bf Proof of Theorem 2.3.}  By Eq. (1.4), one can be convinced that
\begin{eqnarray*}
Kf^{*}(PG_{n+1})&=&\sum_{\{u,v\}\subseteq V_{PG_{n}}}d(u)d(v)r(u,v)+\sum_{v\in V_{PG_{n}}}\sum_{x_{i}\in V_{O_{n+1}}}d(v)d(x_{i})r(v,x_{i})\\
&&+\sum_{\{x_{i},x_{j}\}\subseteq V_{O_{n+1}}}d(x_{i})d(x_{j})r(x_{i},x_{j}).
\end{eqnarray*}
Note that
\begin{align*}
&\sum_{\{u,v\}\subseteq V_{PG_{n}}}d(u)d(v)r(u,v)=Kf^{*}(PG_{n})+\sum_{v\in V_{PG_{n}}}d(v)r(u_{n},v).
\end{align*}
Recall that $d(x_{1})=3$ and $d(x_{i})=2$ for $i \in \{2, ~3,~ 4,~ 5\}$. From Eq. (2.15), we have
\begin{eqnarray*}
\sum_{v\in V_{PG_{n}}}\sum_{x_{i}\in V_{O_{n+1}}}d(v)d(x_{i})r(v,x_{i})&=&\sum_{v\in V_{PG_{n}}}d(v)\bigg[3\big(r(u_{n},v)+1\big)+4\big(r(u_{n},v)+1+\frac{4}{5}\big)\\
&&+4\big(r(u_{n},v)+1+\frac{6}{5}\big)\bigg]\\
&=&\sum_{v\in V_{PG_{n}}}d(v)\big(11r(u_{n},v)+19\big)\\
&=&11\sum_{v\in V_{PG_{n}}}d(v)r(u_{n},v)+19(12n-1).
\end{eqnarray*}
From Eq. (2.16), one follows that
\begin{eqnarray*}
\sum_{\{x_{i},x_{j}\}\subseteq V_{O_{n+1}}}d(x_{i})d(x_{j})r(x_{i},x_{j})&=&\frac{1}{2}\sum^{5}_{i=1}d(x_{i})\big(\sum^{5}_{j=1}d(x_{j})r(x_{i},x_{j})\big)\\
&=&\frac{1}{2}(3\times8+2\times 2\times \frac{44}{5}+2\times 2\times \frac{46}{5})\\
&=&96.
\end{eqnarray*}
Then
\begin{eqnarray}
Kf^{*}(PG_{n+1})=Kf^{*}(PG_{n})+12\sum_{v\in V_{PG_{n}}}d(v)r(u_{n},v)+228n+77.
\end{eqnarray}

For the random chain networks $PG_{n}$, we obtain that $\sum_{v\in V_{PG_{n}}}d(v)r(u_{n},v)$ is a stochastic variate. Let
\begin{eqnarray*}
C_{n}:=\mathbb{E}\big(\sum_{v\in V_{PG_{n}}}d(v)r(u_{n},v)\big).
\end{eqnarray*}

Then, let's continue to consider the two probable cases shown below.

{\bf Case 1.} $PG_{n}\longrightarrow PG^{1}_{n+1}$.

In this case, $u_{n}$ (of $PG_{n}$) overlaps with $x_{2}$ or $x_{5}$ (of $O_{n}$). Therefore, $\sum\limits_{v\in V_{PG_{n}}}d(v)r(u_{n},v)$ is rewritten as $\sum\limits_{v\in V_{PG_{n}}}d(v)r(x_{2},v)$ or $\sum\limits_{v\in V_{PG_{n}}}d(v)r(x_{5},v)$ with probability $p_{1}$.

{\bf Case 2.} $PG_{n}\longrightarrow PG^{2}_{n+1}$.

In this case, $u_{n}$ (of $PG_{n}$) overlaps with $x_{3}$ or $x_{4}$ (of $O_{n}$). Therefore, $\sum\limits_{v\in V_{PG_{n}}}d(v)r(u_{n},v)$ is rewritten as $\sum\limits_{v\in V_{PG_{n}}}d(v)r(x_{3},v)$ or $\sum\limits_{v\in V_{PG_{n}}}d(v)r(x_{4},v)$ with probability $p_{2}=1-P_{1}$.

Together with the above Cases, by applying the expectation operator and Eq. (2.15), Eq. (2.16), one follows that
\begin{eqnarray*}
C_{n}&=&p_{1}\sum\limits_{v\in V_{PG_{n}}}d(v)r(x_{2},v)+(1-p_{1})\sum\limits_{v\in V_{PG_{n}}}d(v)r(x_{3},v)\\
&=&p_{1}\big[\sum\limits_{v\in V_{PG_{n-1}}}d(v)r(u_{n-1},v)+(1+\frac{4}{5})\times(12n-13)+\frac{44}{5}\big]\\
&&+(1-p_{1})\big[\sum\limits_{v\in V_{PG_{n-1}}}d(v)r(u_{n-1},v)+(1+\frac{6}{5})\times(12n-13)+\frac{46}{5}\big].
\end{eqnarray*}
By taking the expectation of both sides of the above equation, and noting that $\mathbb{E}(C_{n})=C_{n}$, we obtain
\begin{eqnarray*}
C_{n}&=&C_{n-1}+p_{1}\big[\frac{9}{5}\times(12n-13)+\frac{44}{5}\big]\\
&&+(1-p_{1})\big[\frac{11}{5}\times(12n-13)+\frac{46}{5}\big]\\
&=&C_{n-1}+(\frac{132}{5}-\frac{24}{5}p_{1})n+(\frac{24}{5}p_{1}-\frac{97}{5}).
\end{eqnarray*}
Meanwhile, for $n=1$, the boundary condition is
\begin{eqnarray*}
C_{1}=\mathbb{E}\big(\sum_{v\in V_{PG_{1}}}d(v)r(u_{1},v)\big)=8.
\end{eqnarray*}
Using above condition and the recurrence relation with respect to $C_{n}$, it is no hard to obtain
\begin{eqnarray}
C_{n}&=&(\frac{66}{5}-\frac{12}{5}p_{1})n^{2}+(\frac{12}{5}p_{1}-\frac{31}{5})n+1.
\end{eqnarray}
So by combining formula (2.17), (2.18), we can obtain
\begin{eqnarray*}
\mathbb{E}\big(Kf^{*}(PG_{n+1})\big)&=&\mathbb{E}\big(Kf^{*}(PG_{n})\big)+12C_{n}+228n+77\\
&=&\mathbb{E}\big(Kf^{*}(PG_{n})\big)+12\big((\frac{66}{5}-\frac{12}{5}p_{1})n^{2}+(\frac{12}{5}p_{1}-\frac{31}{5})n+1\big)+228n+77\\
&=&\mathbb{E}\big(Kf^{*}(PG_{n})\big)+(\frac{792}{5}-\frac{144}{5}p_{1})n^{2}+(\frac{144}{5}p_{1}+\frac{768}{5})n+89.
\end{eqnarray*}
For $n=1$, we obtain $\mathbb{E}\big(Kf^{*}(PG_{1})\big)=40$.

Similarly, according to the recurrence relation related to $\mathbb{E}\big(Kf^{*}(PG_{n})\big)$, we have
\begin{eqnarray*}
\mathbb{E}\big(Kf^{*}(PG_{n})\big)&=&(\frac{264}{5}-\frac{48}{5}p_{1})n^{3}+(\frac{144}{5}p_{1}-\frac{12}{5})n^{2}+(\frac{193}{5}-\frac{96}{5}p_{1})n-49,
\end{eqnarray*}
as desired.\hfill\rule{1ex}{1ex}\

\begin{thm}
For $n\geq 1$, the analytic expression of $\mathbb{E}\big(Kf^{+}(PG_{n})\big)$ for the random chain networks is
\begin{eqnarray*}
\mathbb{E}\big(Kf^{+}(PG_{n})\big)&=&(44-8p_{1})n^{3}+(48p_{1}+11)n^{2}-(88p_{1}+15)n+48p_{1}.
\end{eqnarray*}
\end{thm}

\noindent{\bf Proof of Theorem 2.4.} Notice that the random chain networks $PG_{n+1}$ is organized by adding a new terminal pentagon $O_{n+1}$ to $PG_{n}$ by an edge. By Eq. (1.5), one has
\begin{eqnarray}
Kf^{+}(PG_{n+1})=\Theta_{1}+\Theta_{2}+\Theta_{3},
\end{eqnarray}
where
\begin{eqnarray*}
\Theta_{1}=\sum_{\{u,v\}\subseteq V_{PG_{n}}}\big(d(u)+d(v)\big)r(u,v),\\
\Theta_{2}=\sum_{v\in V_{PG_{n}}}\sum_{x_{i}\in V_{O_{n+1}}}\big(d(v)+d(x_{i})\big)r(v,x_{i}),\\
\Theta_{3}=\sum_{\{x_{i},x_{j}\}\subseteq V_{O_{n+1}}}\big(d(x_{i})+d(x_{j})\big)r(x_{i},x_{j}).\\
\end{eqnarray*}
Let $r(u_{n}\mid PG_{n}):=\sum\limits_{v\in V_{PG_{n}}}r(u_{n},v)$. It is routine to check that
\begin{eqnarray*}
\Theta_{1}&=&\sum_{\{u,v\}\subseteq V_{PG_{n}}}\big(d(u)+d(v)\big)r(u,v)=Kf^{+}(PG_{n})+r(u_{n}\mid PG_{n})
\end{eqnarray*}
Note that $PG_{n}$ has $5n$ vertices. Note that $d(x_{1})=3$, $d(x_{i})=2$ for $i \in \{2,~ 3, ~4, ~5\}$. By Eq. (2.15), we can get
\begin{eqnarray*}
\Theta_{2}&=&\sum_{v\in V_{PG_{n}}}\sum_{x_{i}\in V_{O_{n+1}}}\big(d(v)+d(x_{i})\big)r(v,x_{i})\\
&=&\sum_{v\in V_{PG_{n}}}\sum_{x_{i}\in V_{O_{n+1}}}d(v)r(v,x_{i})+\sum_{v\in V_{PG_{n}}}\sum_{x_{i}\in V_{O_{n+1}}}d(x_{i})r(v,x_{i})\\
&=&\sum_{v\in V_{PG_{n}}}d(v)\bigg[\big(r(u_{n},v)+1\big)+\big(r(u_{n},v)+1+\frac{4}{5}\big)+\big(r(u_{n},v)+1+\frac{6}{5}\big)+\big(r(u_{n},v)+1+\frac{6}{5}\big)\\
&&+\big(r(u_{n},v)+1+\frac{4}{5}\big)\bigg]+\sum_{v\in V_{PG_{n}}}\bigg[3\big(r(u_{n},v)+1\big)+2\big(r(u_{n},v)+1+\frac{4}{5}\big)+2\big(r(u_{n},v)+1+\frac{6}{5}\big)\\
&&+2\big(r(u_{n},v)+1+\frac{6}{5}\big)+2\big(r(u_{n},v)+1+\frac{4}{5}\big)\bigg]\\
&=&5\sum_{v\in V_{PG_{n}}}d(v)r(u_{n},v)+9(12n-1)+11r(u_{n}\mid PG_{n})+19\times5n.
\end{eqnarray*}
Note that $\sum\limits_{i=1}^{5}r(x_{k},x_{i})=4$ for $k=1, ~2,\cdot \cdot \cdot, ~5$. From Eq. (2.16), one sees that
\begin{eqnarray*}
\Theta_{3}&=&\sum_{\{x_{i},x_{j}\}\subseteq V_{O_{n+1}}}\big(d(x_{i})+d(x_{j})\big)r(x_{i},x_{j})=\frac{1}{2}\sum^{5}_{i=1}\sum^{5}_{j=1}\big(d(x_{i})+d(x_{j})\big)r(x_{i},x_{j})\\
&=&\sum^{5}_{i=1}\sum^{5}_{j=1}d(x_{i})r(x_{i},x_{j})=4\times(3+2\times 4)=44.
\end{eqnarray*}
Then, Eq. (2.19) can be rewritten as
\begin{eqnarray}
Kf^{+}(PG_{n+1})=Kf^{+}(PG_{n})+5\sum_{v\in V_{PG_{n}}}d(v)r(u_{n},v)+12r(u_{n}\mid PG_{n})+203n+35.
\end{eqnarray}

For a the random chain networks $PG_{n}$, we know that $r(u_{n}\mid PG_{n})$ is a random variable. Let
\begin{eqnarray*}
D_{n}:=\mathbb{E}\big(r(u_{n}\mid PG_{n})\big).
\end{eqnarray*}

We proceed by taking into account the following two cases.

{\bf Case 1.} $PG_{n}\longrightarrow PG^{1}_{n+1}$.

In this case, $u_{n}$ (of $PG_{n}$) coincides with $x_{2}$ or $x_{5}$ (of $O_{n}$). Then, $r(u_{n}\mid PG_{n})$ is given by $r(x_{2}\mid PG_{n})$ or $r(x_{5}\mid PG_{n})$ with probability $p_{1}$.

{\bf Case 2.} $PG_{n}\longrightarrow PG^{2}_{n+1}$.

In this case, $u_{n}$ (of $PG_{n}$) coincides with $x_{3}$ or $x_{4}$ (of $O_{n}$). Then, $r(u_{n}\mid PG_{n})$ is given by $r(x_{3}\mid PG_{n})$ or $r(x_{4}\mid PG_{n})$ with probability $p_{2}=1-p_{1}$.

Together with Cases 1 and 2, we obtain
\begin{eqnarray*}
D_{n}&=&p_{1}\cdot r(x_{2}\mid PG_{n})+(1-p_{1})\cdot r(x_{3}\mid PG_{n})\\
&=&p_{1}\cdot \big[r(u_{n-1}\mid PG_{n-1})+5(n-1)(1+\frac{4}{5})+4\big]+(1-p_{1})\cdot \big[r(u_{n-1}\mid PG_{n-1})\\
&&+5(n-1)(1+\frac{6}{5})+4\big].
\end{eqnarray*}
By taking the expectation of both sides of the above equation, and noting that $\mathbb{E}(D_{n})=D_{n}$, we obtain
\begin{eqnarray*}
D_{n}&=&p_{1}\big[D_{n-1}+9(n-1)+4\big]+(1-p_{1})\big[D_{n-1}+11(n-1)+4\big]\\
&=&D_{n-1}+(11-2p_{1})n+(6p_{1}-7).
\end{eqnarray*}
Meanwhile, for $n=1$, the boundary condition is
\begin{eqnarray*}
D_{1}=\mathbb{E}\big(r(u_{1}\mid PG_{1})\big)=4.
\end{eqnarray*}
Using above condition and the recurrence relation with respect to $D_{n}$, it is no hard to obtain
\begin{eqnarray}
D_{n}&=&(\frac{11}{2}-p_{1})n^{2}+(5p_{1}-\frac{3}{2})n-4p_{1}.
\end{eqnarray}
So by combining formula (2.20), (2.21) we can obtain
\begin{eqnarray*}
\mathbb{E}\big(Kf^{+}(PG_{n})\big)&=&\mathbb{E}\big(Kf^{+}(PG_{n})\big)+12D_{n}+5C_{n}+203n+35\\
&=&\mathbb{E}\big(Kf^{+}(PG_{n})\big)+12\big[(\frac{11}{2}-p_{1})n^{2}+(5p_{1}-\frac{3}{2})n-4p_{1}\big]+5\big[(\frac{66}{5}-\frac{12}{5}p_{1})n^{2}\\
&&+(\frac{12}{5}p_{1}-\frac{31}{5})n+1\big]+203n+35\\
&=&\mathbb{E}\big(Kf^{+}(PG_{n})\big)+(132-24p_{1})n^{2}+(72p_{1}+154)n+(40-48p_{1}).
\end{eqnarray*}
For $n=1$, we obtain $\mathbb{E}\big(Kf^{+}(PG_{1})\big)=40$.

Similarly, according to the recurrence relation related to $\mathbb{E}\big(Kf^{+}(PG_{n})\big)$, we have
\begin{eqnarray*}
\mathbb{E}\big(Kf^{+}(PG_{n})\big)&=&(44-8p_{1})n^{3}+(48p_{1}+11)n^{2}-(88p_{1}+15)n+48p_{1},
\end{eqnarray*}
as desired.\hfill\rule{1ex}{1ex}\

Next, we use the computer to represent the images of expectation functions of the four indices about the random chain networks. Through the image, we can see that the expected values of the four topological indices of the random chain networks are related to two variables $p_{1}$ and $n$, and all of them show a positive correlation. Among them, the growth rate of $\mathbb{E}\big(Gut(PG_{n})\big)$ is the fastest and that of $\mathbb{E}\big(Kf^{*}(PG_{n})\big)$ is the slowest, while that of $\mathbb{E}\big(S(PG_{n})\big)$ and $\mathbb{E}\big(Kf^{+}(PG_{n})\big)$ is between them.

\begin{figure}[htbp]
\centering\includegraphics[width=12.936cm,height=9.548cm]{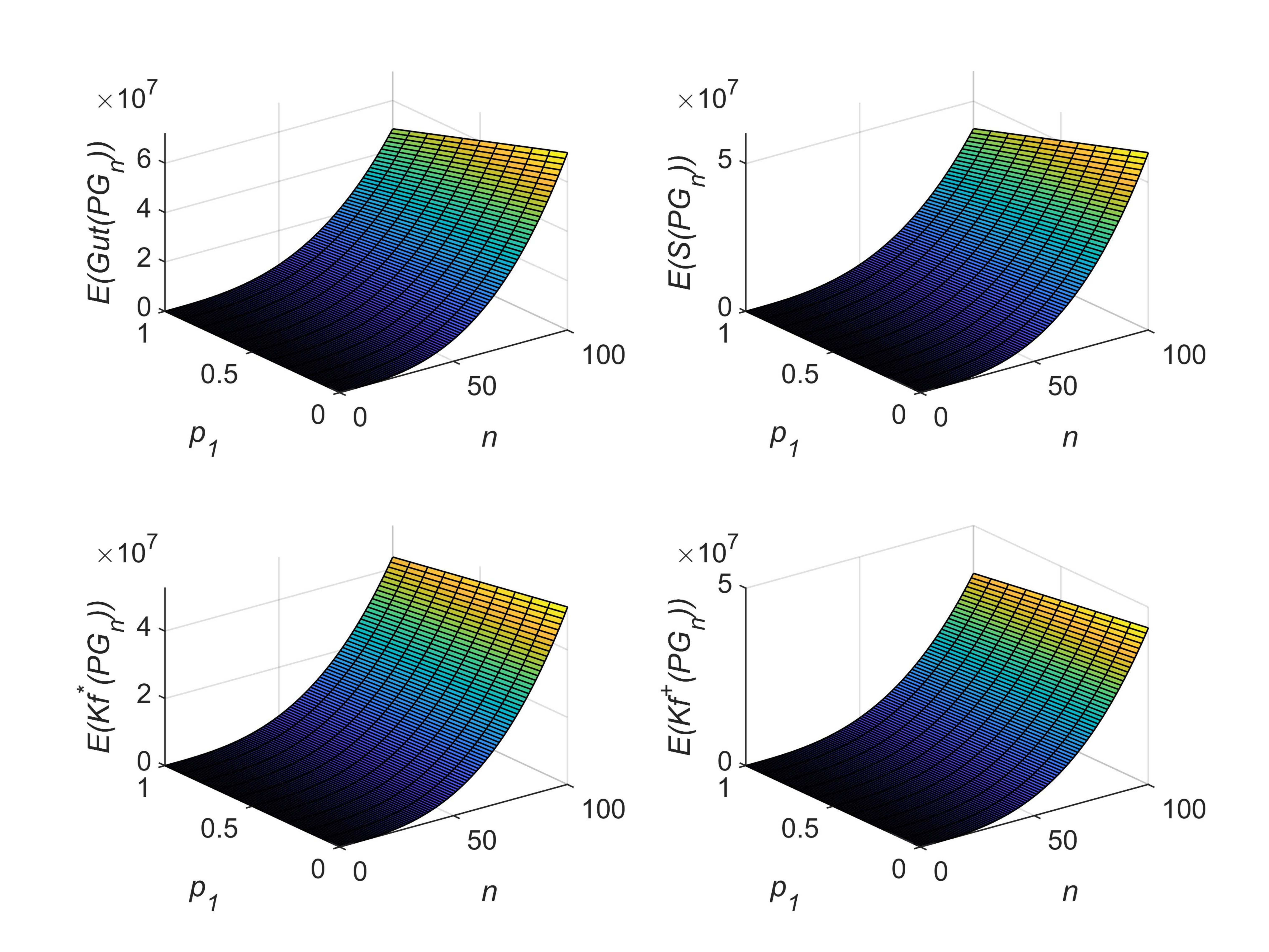}
\caption{ The expected function images of the four indices.}
\end{figure}

\section{The limiting behaviours for $Gut(PG_{n})$, $S(PG_{n})$, $Kf^{*}(PG_{n})$ and $Kf^{+}(PG_{n})$ of the random chain networks}
In Section 3, we obtain the explicit analytical expressions for $Var\big(Gut(PG_{n})\big)$, $Var\big(S(PG_{n})\big)$, $Var\big(Kf^{*}(PG_{n})\big)$ and $Var\big(Kf^{+}(PG_{n})\big)$. For the random chain networks, we prove the Gutman, Schultz, multiplicative degree-Kirchhoff and additive degree-kirchhoff indices asymptotically obey normal distributions. We use the same notation as those used at Section 2.

\begin{thm}
Suppose Hypotheses 1 and 2 are correct, then the following results are also right.

$(i)$ The variance of the Gutman index is denoted by
\begin{eqnarray*}
Var\big(Gut(PG_{n})\big)&=&\frac{1}{30}\big({\sigma}^{2}_{1}n^{5}-5r_{1}n^{4}+10\tilde{\sigma}^{2}_{1}n^{3}+(65r_{1}-30{\sigma}^{2}_{1}-45\tilde{\sigma}^{2}_{1})n^{2}\\
&&+(59{\sigma}^{2}_{1}+65\tilde{\sigma}^{2}_{1}-120r_{1})n+(60r_{1}-30{\sigma}^{2}_{1}-30\tilde{\sigma}^{2}_{1})\big),
\end{eqnarray*}
where
\begin{align*}
\sigma^{2}_{1}&=288^{2}p_{1}+432^{2}(1-p_{1})-\big(288p_{1}+432(1-p_{1})\big)^{2},\\
\tilde{\sigma}^{2}_{1}&=156^{2}p_{1}+300^{2}(1-p_{1})-\big(156p_{1}+300(1-p_{1})\big)^{2},\\
r_{1}&=288\cdot156\cdot p_{1}+432\cdot300\cdot(1-p_{1})-\big(288p_{1}+432(1-p_{1})\big)\cdot\big(156p_{1}+300(1-p_{1})\big).
\end{align*}

$(ii)$ For $n\rightarrow \infty$, $Gut(PG_{n})$ asymptotically obeys normal distributions. One has
\begin{eqnarray*}
\lim_{n\rightarrow \infty}\sup_{a\in \mathbb{R}}\mid \mathbb{P}\bigg(\frac{Gut(PG_{n})-\mathbb{E}\big(Gut(PG_{n})\big)}{\sqrt{Var\big(Gut(PG_{n})\big)}}\leq a\bigg)-{\int}^{a}_{-\infty}\frac{1}{\sqrt{2\pi}}e^{-\frac{t^{2}}{2}}dt \mid =0.
\end{eqnarray*}
\end{thm}

\noindent{\bf Proof of Theorem 3.1.} Let $E_{n}=12\sum\limits_{v \in V_{PG_{n}}}d(v)d(u_{n},v)$. Then by Eq. (2.9), we obtain
\begin{eqnarray}
Gut(PG_{n+1})=Gut(PG_{n})+E_{n}+276n+49.
\end{eqnarray}

Recalling that $Z_{n}^{1}$ and $Z_{n}^{2}$ are random variables which stand for our choice to construct $PG_{n+1}$ by $PG_{n}$. We have the next two facts.

\noindent{\bf Fact 3.1.1.} $E_{n}Z_{n}^{1}=(E_{n-1}+288n-156)Z_{n}^{1}.$

\noindent{\bf Proof.} If $Z_{n}^{1}$ = 0, the result is obvious. Then, we only take into account $Z_{n}^{1}$ = 1, which implies $PG_{n} \rightarrow PG^{1}_{n+1}$. In this case, $u_{n}$ (of $PG_{n}$) overlaps with $x_{2}$ or $x_{5}$ (of $O_{n}$), see Fig. 3. In this situation, by using Eq. (2.7) and Eq. (2.8),
\begin{eqnarray*}
E_{n}&=&12\sum_{v \in V_{PG_{n}}}d(v)d(x_{2},v)\\
&=&12\sum_{v \in V_{PG_{n-1}}}d(v)d(x_{2},v)+12\sum_{v \in V_{O_{n}}}d(v)d(x_{2},v)\\
&=&12\sum_{v \in V_{PG_{n-1}}}d(v)\big(d(v,u_{n-1})+d(x_{2},u_{n-1})\big)+12\times 13\\
&=&12\sum_{v \in V_{PG_{n-1}}}d(v)\big(d(v,u_{n-1})+2\big)+12\times 13\\
&=&E_{n-1}+24\times(12n-13)+156\\
&=&E_{n-1}+288n-156.
\end{eqnarray*}
Thus, we conclude the desired Fact.

\noindent{\bf Fact 3.1.2.} $E_{n}Z_{n}^{2}=(E_{n-1}+432n-300)Z_{n}^{2}.$

Similar to the proof of Fact 3.1.1, we only consider the fact $Z_{n}^{2}=1$, that is $PG_{n} \rightarrow PG^{2}_{n+1}$. In the same way, we omit the details.

Noting that $Z_{n}^{1}+Z_{n}^{2}=1$, by the above discussions, it holds that
\begin{eqnarray*}
E_{n}&=&E_{n}(Z_{n}^{1}+Z_{n}^{2})\\
&=&(E_{n-1}+288n-156)Z_{n}^{1}+(E_{n-1}+432n-300)Z_{n}^{2}\\
&=&E_{n-1}+(288Z_{n}^{1}+432Z_{n}^{2})n-(156Z_{n}^{1}+300Z_{n}^{2})\\
&=&E_{n-1}+nU_{n}-V_{n},
\end{eqnarray*}
where for each $n$,

~~~~~~~~~$U_{n}=288Z_{n}^{1}+432Z_{n}^{2}$,~~~~$V_{n}=156Z_{n}^{1}+300Z_{n}^{2}$.
\\\text{Therefore, by Eq. (3.22), it follows that}\par

\begin{align}
Gut(PG_{n})=& Gut(PG_{1})+\sum_{l=1}^{n-1}E_{l}+\sum_{l=1}^{n-1}(276l+49)\nonumber\\
=& Gut(PG_{1})+\sum_{l=1}^{n-1}(\sum_{m=1}^{l-1}(E_{m+1}-E_{m})+E_{1})+\sum_{l=1}^{n-1}(276l+49)\nonumber\\
=& Gut(PG_{1})+\sum_{l=1}^{n-1}\sum_{m=1}^{l-1}(E_{m+1}-E_{m})+(n-1)E_{1}+\sum_{l=1}^{n-1}(276l+49)\nonumber\\
=& Gut(PG_{1})+\sum_{l=1}^{n-1}\sum_{m=1}^{l-1}\big((m+1)U_{m+1}-V_{m+1}\big)+O(n^2).
\end{align}
By direct calculation, we put

~~~~~~~~~~$Var(U_{m})=\sigma^{2}_{1}$,~~~~$Var(V_{m})=\tilde{\sigma}^{2}_{1}$,~~~~$Cov(U_{m},V_{m})=r_{1}$,
\\\text{where for any two stochastic variates $X$ and $Y$, $Cov(X,Y)=\mathbb{E}(XY)-\mathbb{E}(X)\mathbb{E}(Y)$.}\par

Refer to ref.\cite{45}, by the properties of variance, Eq. (3.23) and exchanging the order of $l$ and $m$, we can directly find out that
\begin{align*}
&Var\big(Gut(PG_{n})\big)\\
&=Var\bigg[\sum_{l=1}^{n-1}\sum_{m=1}^{l-1}\big((m+1)U_{m+1}-V_{m+1}\big)\bigg]=Var\bigg[\sum_{m=1}^{n-2}\sum_{l=m+1}^{n-1}\big((m+1)U_{m+1}-V_{m+1}\big)\bigg]\\
&=Var\bigg[\sum_{m=1}^{n-2}\big((m+1)U_{m+1}-V_{m+1}\big)(n-m-1)\big]=\sum_{m=1}^{n-2}(n-m-1)^{2}Var\big((m+1)U_{m+1}-V_{m+1}\big)\\
&=\sum_{m=1}^{n-2}(n-m-1)^{2}Cov\big((m+1)U_{m+1}-V_{m+1},(m+1)U_{m+1}-V_{m+1}\big)\\
&=\sum_{m=1}^{n-2}(n-m-1)^{2}\big((m+1)^2Cov(U_{m+1},U_{m+1})-2(m+1)Cov(U_{m+1},V_{m+1})+Cov(V_{m+1},V_{m+1})\big)\\
&=\sum_{m=1}^{n-2}(n-m-1)^{2}\big((m+1)^2\sigma^{2}_{1}-2(m+1)r_{1}+\tilde{\sigma}^{2}_{1}\big).\\
\end{align*}
By using a computer, the above expression indicates the result $Theorem ~3.1.~~ (i)$.

Now we will go on to the proof of $Theorem ~3.1.~~ (ii)$. Firstly, for any $n\in \mathbb{N}$, let

~~~~~$\mathcal{U}_{n}=\sum\limits_{l=1}^{n-1}\sum\limits_{m=1}^{l-1}(m+1)U_{b+1}$,~~~~$\mathcal{V}_{n}=\sum\limits_{l=1}^{n-1}\sum\limits_{m=1}^{l-1}V_{m+1}$, ~~$\mu= \mathbb{E}(U_{m})$, ~~and~~   $\phi (t)=\mathbb{E}(e^{t(U_{m}-\mu)})$.
\\\text{By these notations, obviously, we have}\par

$e^{t\big(\mathcal{U}_{n}-\mathbb{E}(\mathcal{U}_{n})\big)}=e^{t\sum\limits_{l=1}^{n-1}\sum\limits_{m=1}^{l-1}(m+1)(U_{m+1}-\mu)}=e^{t\sum\limits_{m=1}^{n-2}\sum\limits_{l=m+1}^{n-1}(m+1)(U_{m+1}-\mu)}=e^{t\sum\limits_{m=1}^{n-2}(n-m-1)(m+1)(U_{m+1}-\mu)}$,
\\\text{then}\par
\begin{align}
\mathbb{E}\bigg[e^{t\big(\mathcal{U}_{n}-\mathbb{E}(\mathcal{U}_{n})\big)}\bigg]=&\mathbb{E}\big(e^{t\sum\limits_{m=1}^{n-2}(n-m-1)(m+1)(U_{m+1}-\mu)}\big)=\prod\limits_{m=1}^{n-2}\mathbb{E}\big(e^{t(n-m-1)(m+1)(U_{m+1}-\mu)}\big)\nonumber\\
=&\prod\limits_{m=1}^{n-2}\phi \big(t(n-m-1)(m+1)\big),
\end{align}
and for some $k>0$,
\begin{eqnarray}
\mathcal{V}_{n}=\sum\limits_{l=1}^{n-1}\sum\limits_{m=1}^{l-1}V_{m+1}\leq kn^2.
\end{eqnarray}
Noting that

~~$Var\big(Gut(PG_{n})\big) \asymp \frac{1}{30}\sigma_{1}^2n^5$,~~ $\phi (t)=1+\frac{\sigma_{1}^2}{2}t^2+O(t^2)$, ~~and ~~$\sum\limits^{n-2}_{m=1}(m+1)^2(n-m-1)^2\asymp \frac{n^5}{30}$.
\\\text{By Taylor's formula and Eqs. (3.23)-(3.25), one holds that}\par
\begin{align*}
&\lim\limits_{n\rightarrow \infty}\mathbb{E}\exp\bigg\{t\frac{Gut(PG_{n})-\mathbb{E}\big(Gut(PG_{n})\big)}{\sqrt{Var\big(Gut(PG_{n})\big)}}\bigg\}\\
&=\lim\limits_{n\rightarrow \infty}\mathbb{E}\exp\bigg\{t\frac{\big(Gut(PG_{1})+\mathcal{U}_{n}-\mathcal{V}_{n}+O(n^2)\big)-\mathbb{E}\big(Gut(PG_{1})+\mathcal{U}_{n}-\mathcal{V}_{n}+O(n^2)\big)}{\frac{\sigma_{1} n^{\frac{5}{2}}}{\sqrt{30}}}\bigg\}\\
&=\lim\limits_{n\rightarrow \infty}\mathbb{E}\exp\bigg\{t\frac{\sqrt{30}\big(\mathcal{U}_{n}-\mathbb{E}(\mathcal{U}_{n})\big)}{\sigma_{1} n^{\frac{5}{2}}}\bigg\}=\lim\limits_{n\rightarrow \infty}\prod\limits_{m=1}^{n}\phi \big(\frac{\sqrt{30}t(m+1)(n-m-1)}{\sigma_{1} n^{\frac{5}{2}}}\big)\\
&=\lim\limits_{n\rightarrow \infty}\exp\bigg\{\sum\limits_{m=1}^{n}\ln\phi \big(\frac{\sqrt{30}t(m+1)(n-m-1)}{\sigma_{1} n^{\frac{5}{2}}}\big)\bigg\}\\
&=\lim\limits_{n\rightarrow \infty}\exp\bigg\{\sum\limits_{m=1}^{n}\ln\big(1+\frac{\sigma_{1}^2}{2}\cdot\frac{30t^2(m+1)^2(n-m-1)^2}{\sigma_{1}^2n^5}+O(\frac{1}{n})\big)\bigg\}\\
&=\lim\limits_{n\rightarrow \infty}\exp\bigg\{\sum\limits_{m=1}^{n-2}\big(\frac{\sigma_{1}^2}{2}\cdot\frac{30t^2(m+1)^2(n-m-1)^2}{\sigma_{1}^2n^5}+O(\frac{1}{n})\big)\bigg\}\\
&=e^{\frac{t^2}{2}}.
\end{align*}
 Assume that $\mathbb{I}$ is a complex number with $\mathbb{I}^2=-1$. We use $\mathbb{I}t$ instead of $t$ and one has
\begin{align*}
\lim\limits_{n\rightarrow \infty}\mathbb{E}\exp\bigg\{\mathbb{I}t\frac{Gut(PG_{n})-\mathbb{E}\big(Gut(PG_{n})\big)}{\sqrt{Var\big(Gut(PG_{n})\big)}}\bigg\}=e^{-\frac{t^2}{2}}.
\end{align*}

According to the above formula (\cite{49}, Chapter 1) and the theory of continuity of probability characteristic functions (\cite{50}, Chapter 15), we complete the proof of $Theorem ~3.1.~~ (ii)$.\hfill\rule{1ex}{1ex}\

\begin{thm}
Suppose Hypotheses 1 and 2 are true, then the next results hold.

 $(i)$ The variance of the Schultz index is denoted by
\begin{eqnarray*}
Var\big(S(PG_{n})\big)&=&\frac{1}{30}\big({\sigma}_{2}^{2}n^{5}-5r_{2}n^{4}+10\tilde{\sigma}^{2}_{2}n^{3}+(65r_{2}-30{\sigma}^{2}_{2}-45\tilde{\sigma}^{2}_{2})n^{2}\\
&&+(59{\sigma}^{2}_{2}+65\tilde{\sigma}^{2}_{2}-120r_{2})n+(60r_{2}-30{\sigma}^{2}_{2}-30\tilde{\sigma}^{2}_{2})\big),
\end{eqnarray*}
where
\begin{align*}
\sigma^{2}_{2}&=240^{2}p_{1}+360^{2}(1-p_{1})-\big(240p_{1}+360(1-p_{1})\big)^{2},\\
\tilde{\sigma}_{2}^{2}&=113^{2}p_{1}+233^{2}(1-p_{1})-\big(113p_{1}+233(1-p_{1})\big)^{2},\\
r_{2}&=240\cdot113p_{1}+360\cdot233(1-p_{1})-\big(240p_{1}+360(1-p_{1})\big)\big(113p_{1}+233(1-p_{1})\big).
\end{align*}

$(ii)$  For $n\rightarrow \infty$, $S(PG_{n})$ asymptotically obeys normal distributions. One has
\begin{eqnarray*}
\lim_{n\rightarrow \infty}\sup_{a\in \mathbb{R}}\mid \mathbb{P}\bigg(\frac{S(PG_{n})-\mathbb{E}\big(S(PG_{n})\big)}{\sqrt{Var\big(S(PG_{n})\big)}}\leq a\bigg)-{\int}^{a}_{-\infty}\frac{1}{\sqrt{2\pi}}e^{-\frac{t^{2}}{2}}dt \mid =0.
\end{eqnarray*}
\end{thm}

\noindent{\bf Proof of Theorem 3.2.} Let $F_{n}=\sum\limits_{v \in V_{PG_{n}}}\big(5d(v)+12\big)d(u_{n},v)$. Then by Eq. (2.13), we obtain
\begin{eqnarray}
S(PG_{n+1})=S(PG_{n})+F_{n}+247n+55.
\end{eqnarray}

According to the previous proof of Theorem 3.1, the two facts are obtained.

\noindent{\bf Fact 3.2.1.} $F_{n}Z_{n}^{1}=(F_{n-1}+240n-113)Z_{n}^{1}.$

\noindent{\bf Proof.} If $Z_{n}^{1}$ = 0, the above result is distinct. So we take into account $Z_{n}^{1}$ = 1, which indicates $PG_{n} \rightarrow PG^{1}_{n+1}$. In this case, $u_{n}$ (of $PG_{n}$) coincides with $x_{2}$ or $x_{5}$ (of $O_{n}$), see Fig. 3. In this situation, by using Eq. (2.7) and Eq. (2.8), $F_{n}$ becomes
\begin{align*}
&\sum\limits_{v \in V_{PG_{n}}}\big(5d(v)+12\big)d(x_{2},v)\\
&=\sum_{v \in V_{PG_{n-1}}}\big(5d(v)+12\big)d(x_{2},v)+\sum_{v \in V_{O_{n}}}\big(5d(v)+12\big)d(x_{2},v)\\
&=\sum_{v \in V_{PG_{n-1}}}\big(5d(v)+12\big)\big(d(u_{n-1},v)+2\big)+5\times 13+12\times 6\\
&=\sum_{v \in V_{PG_{n-1}}}\big(5d(v)+12\big)d(u_{n-1},v)+2\cdot \big(5\times12(n-13)+12\times5(n-1)\big)+137\\
&=F_{n-1}+240n-113.
\end{align*}
Thus, we conclude the desired Fact.

\noindent{\bf Fact 3.2.2.} $F_{n}Z_{n}^{2}=(F_{n-1}+360n-233)Z_{n}^{2}.$

As that in the proof of Fact 3.2.1, we only consider the fact $Z_{n}^{2}=1$, that is $PG_{n} \rightarrow PG^{2}_{n+1}$. The proof is similar and details are omitted.

Noting that $Z_{n}^{1}+Z_{n}^{2}=1$, by the above discussions, it holds that
\begin{eqnarray*}
F_{n}&=&F_{n}(Z_{n}^{1}+Z_{n}^{2})\\
&=&(F_{n-1}+240n-113)Z_{n}^{1}+(F_{n-1}+360n-233)Z_{n}^{2}\\
&=&F_{n-1}+(240Z_{n}^{1}+360Z_{n}^{2})n-(113Z_{n}^{1}+233Z_{n}^{2})\\
&=&F_{n-1}+nU^{'}_{n}-V^{'}_{n},
\end{eqnarray*}
where for each $n$,

~~~~~~~~~$U^{'}_{n}=240Z_{n}^{1}+360Z_{n}^{2}$,~~~~$V^{'}_{n}=113Z_{n}^{1}+233Z_{n}^{2}$.
\\\text{Therefore, by Eq. (3.26),}\par

\begin{align}
S(PG_{n})=& S(PG_{1})+\sum_{l=1}^{n-1}F_{l}+\sum_{l=1}^{n-1}(247l+55)\nonumber\\
=& S(PG_{1})+\sum_{l=1}^{n-1}\big(\sum_{m=1}^{l-1}(F_{m+1}-F_{m})+F_{1}\big)+\sum_{l=1}^{n-1}(247l+55)\nonumber\\
=& S(PG_{1})+\sum_{l=1}^{n-1}\sum_{m=1}^{l-1}(F_{m+1}-F_{m})+(n-1)F_{1}+\sum_{l=1}^{n-1}(247l+55)\nonumber\\
=& S(PG_{1})+\sum_{l=1}^{n-1}\sum_{m=1}^{l-1}\big((m+1)U^{'}_{m+1}-V^{'}_{m+1}\big)+O(n^2)\nonumber.
\end{align}
Suppose that

~~~~~~~~~$Var(U^{'}_{m})=\sigma^{2}_{2}$,~~~~$Var(V^{'}_{m})=\tilde{\sigma}^{2}_{2}$,~~~~$Cov(U^{'}_{m},V^{'}_{m})=r_{2}$.

If we substitute $Gut(PG_{n})$ by $S(PG_{n})$ in the proof of Theorem 3.1, the rest of the proof is identical to the proof of Theorem 3.1 and the details are omitted here.\hfill\rule{1ex}{1ex}\

\begin{thm}
Suppose Hypotheses 1 and 2 are true, there are the next main results.

$(i)$ The variance of $Kf^{*}(PG_{n})$ is denoted by
\begin{eqnarray*}
Var\big(Kf^{*}(PG_{n})\big)&=&\frac{1}{30}\big({\sigma}^{2}_{3}n^{5}-5r_{3}n^{4}+10\tilde{\sigma}^{2}_{3}n^{3}+(65r_{3}-30{\sigma}^{2}_{3}-45\tilde{\sigma}^{2}_{3})n^{2}\\
&&+(59{\sigma}^{2}_{3}+65\tilde{\sigma}^{2}_{3}-120r_{3})n+(60r_{3}-30{\sigma}^{2}_{3}-30\tilde{\sigma}^{2}_{3})\big),
\end{eqnarray*}
where
\begin{align*}
\sigma^{2}_{3}&=(\frac{1296}{5})^{2}p_{1}+(\frac{1584}{5})^{2}(1-p_{1})-\big(\frac{1296}{5}p_{1}+\frac{1584}{5}(1-p_{1})\big)^{2},\\
\tilde{\sigma}^{2}_{3}&=(\frac{876}{5})^{2}p_{1}+(\frac{1164}{5})^{2}(1-p_{1})-\big(\frac{876}{5}p_{1}+\frac{1164}{5}(1-p_{1})\big)^{2},\\
r_{3}&=\frac{1296}{5}\cdot\frac{876}{5}p_{1}+\frac{1584}{5}\cdot\frac{1164}{5}(1-p_{1})-\big(\frac{1296}{5}p_{1}+\frac{1584}{5}(1-p_{1})\big)\big(\frac{876}{5}p_{1}+\frac{1164}{5}(1-p_{1})\big).
\end{align*}

$(ii)$ For $n\rightarrow \infty$, $Kf^{*}(PG_{n})$ asymptotically obeys normal distributions. One has
\begin{eqnarray*}
\lim_{n\rightarrow \infty}\sup_{a\in \mathbb{R}}\mid \mathbb{P}\bigg(\frac{Kf^{*}(PG_{n})-\mathbb{E}\big(Kf^{*}(PG_{n})\big)}{\sqrt{Var\big(Kf^{*}(PG_{n})\big)}}\leq a\bigg)-{\int}^{a}_{-\infty}\frac{1}{\sqrt{2\pi}}e^{-\frac{t^{2}}{2}}dt \mid =0.
\end{eqnarray*}
\end{thm}

\noindent{\bf Proof of Theorem 3.3.} Obtained by formula (2.17),  we see
\begin{eqnarray*}
Kf^{*}(PG_{n+1})=Kf^{*}(PG_{n})+12\sum_{v\in V_{PG_{n}}}d(v)r(u_{n},v)+288n+77.
\end{eqnarray*}
Let $H_{n}=12\sum\limits_{v \in V_{PG_{n}}}d(v)r(u_{n},v)$, we have
\begin{eqnarray}
Kf^{*}(PG_{n+1})=Kf^{*}(PG_{n})+H_{n}+288n+77.
\end{eqnarray}

Recalling that $Z_{n}^{1}$, and $Z_{n}^{4}$ are random variables which show the way to construct $PG_{n+1}$ from $PG_{n}$. We get the following four Facts.

\noindent{\bf Fact 3.3.1.} $H_{n}Z_{n}^{1}=(H_{n-1}+\frac{1296}{5}n-\frac{876}{5})Z_{n}^{1}.$

\noindent{\bf Proof.} If $Z_{n}^{1}$ = 0, the result is obvious. So we just take into account $Z_{n}^{1}$ = 1, which implies $PG_{n} \rightarrow PG^{1}_{n+1}$. In this case, $u_{n}$ (of $PG_{n}$) overlaps with $x_{2}$ or $x_{5}$ (of $O_{n}$), see Fig. 3. In this situation, by using Eq. (2.15) and Eq. (2.16), $H_{n}$ becomes
\begin{eqnarray*}
12\sum_{v \in V_{PG_{n}}}d(v)r(x_{2},v)&=&12\sum_{v \in V_{PG_{n-1}}}d(v)r(x_{2},v)+12\sum_{v \in V_{O_{n}}}d(v)r(x_{2},v)\\
&=&12\sum_{v \in V_{PG_{n-1}}}d(v)\big(r(u_{n-1},v)+1+\frac{4}{5}\big)+12\times \frac{44}{5}\\
&=&H_{n-1}+12\times \frac{9}{5}\times (12n-13)+12\times \frac{44}{5}\\
&=&H_{n-1}+\frac{1296}{5}n-\frac{876}{5}.
\end{eqnarray*}
Thus, we conclude the desired Fact.

\noindent{\bf Fact 3.3.2.} $H_{n}Z_{n}^{2}=(H_{n-1}+\frac{1584}{5}n-\frac{1164}{5})Z_{n}^{2}.$

As that in the proof of Fact 4.1.1, we only consider the fact $Z_{n}^{2}=1$, that is $PG_{n} \rightarrow PG^{2}_{n+1}$. The proof is similar and details are omitted.

Noting that $Z_{n}^{1}+Z_{n}^{2}=1$, by the above discussions, it holds that
\begin{eqnarray*}
H_{n}&=&H_{n}(Z_{n}^{1}+Z_{n}^{2})\\
&=&(H_{n-1}+\frac{1296}{5}n-\frac{876}{5})Z_{n}^{1}+(P_{n-1}+\frac{1584}{5}n-\frac{1164}{5})Z_{n}^{2}\\
&=&H_{n-1}+(\frac{1296}{5}Z_{n}^{1}+\frac{1584}{5}Z_{n}^{2})n-(\frac{876}{5}Z_{n}^{1}+\frac{1164}{5}Z_{n}^{2})\\
&=&H_{n-1}+n\hat{U}_{n}-\hat{V}_{n},
\end{eqnarray*}
where for each $n$,

~~~~~~~~~$\hat{U}_{n}=\frac{1296}{5}Z_{n}^{1}+\frac{1584}{5}Z_{n}^{2}$,~~~~$\hat{V}_{n}=\frac{876}{5}Z_{n}^{1}+\frac{1164}{5}Z_{n}^{2}$.
\\\text{Therefore, by Eq. (3.27),}\par
\begin{align}
Kf^{*}(PG_{n})=& Kf^{*}(PG_{1})+\sum_{l=1}^{n-1}H_{l}+\sum_{l=1}^{n-1}(228l+77)\nonumber\\
=& Kf^{*}(PG_{1})+\sum_{l=1}^{n-1}\big(\sum_{m=1}^{l-1}(H_{m+1}-H_{m})+H_{1}\big)+\sum_{l=1}^{n-1}(228l+77)\nonumber\\
=& Kf^{*}(PG_{1})+\sum_{l=1}^{n-1}\sum_{m=1}^{l-1}\big((m+1)\hat{U}_{m+1}-\hat{V}_{m+1}\big)+(n-1)P_{1}+\sum_{l=1}^{n-1}(228l+77)\nonumber\\
=& Kf^{*}(PG_{1})+\sum_{l=1}^{n-1}\sum_{m=1}^{l-1}\big((m+1)\hat{U}_{m+1}-\hat{V}_{m+1}\big)+O(n^2).
\end{align}

By direct calculation, we put

~~~~~~~~~$Var(\hat{U}_{m})=\sigma^{2}_{3}$,~~~~$Var(\hat{V}_{m})=\tilde{\sigma}^{2}_{3}$,~~~~$Cov(\hat{U}_{m},\hat{V}_{m})=r_{3}$.

If we substitute $Gut(PG_{n})$ by $S(PG_{n})$ in the proof of Theorem 3.1, the rest of the proof is identical to the proof of Theorem 3.1 and the details are omitted here.\hfill\rule{1ex}{1ex}\

We proceed by showing the following result about the expatiatory formula of the variance of $Kf^{+}(PG_{n})$.

\begin{thm}
 Suppose Hypotheses 1 and 2 are true, there are the next main results.

 $(i)$ The variance of $Kf^{+}(PG_{n})$ is denoted by
\begin{eqnarray*}
Var\big(Kf^{+}(PG_{n})\big)&=&\frac{1}{30}\big({\sigma}^{2}_{4}n^{5}-5r_{4}n^{4}+10\tilde{\sigma}^{2}_{4}n^{3}+(65r_{4}-30{\sigma}^{2}_{4}-45\tilde{\sigma}^{2}_{4})n^{2}\\
&&+(59{\sigma}^{2}_{4}+65\tilde{\sigma}^{2}_{4}-120r_{4})n+(60r_{4}-30{\sigma}^{2}_{4}-30\tilde{\sigma}^{2}_{4})\big),
\end{eqnarray*}
where
\begin{align*}
\sigma^{2}_{4}&=216^{2}p_{1}+264^{2}(1-p_{1})-\big(216p_{1}+264(1-p_{1})\big)^{2},\\
\tilde{\sigma}^{2}_{4}&=133^{2}p_{1}+181^{2}(1-p_{1})-\big(133p_{1}+181(1-p_{1})\big)^{2},\\
r_{4}&=216\cdot133P_{1}+264\cdot181(1-P_{1})-\big(216p_{1}+264(1-p_{1})\big)\big(133p_{1}+181(1-p_{1})\big).
\end{align*}

$(ii)$ For $n \rightarrow \infty$, $Kf^{+}(PG_{n})$ asymptotically obeys normal distributions. One has
\begin{eqnarray*}
\lim_{n\rightarrow \infty}\sup_{a\in \mathbb{R}}\mid \mathbb{P}\bigg(\frac{Kf^{+}(PG_{n})-\mathbb{E}\big(Kf^{+}(PG_{n})\big)}{\sqrt{Var\big(Kf^{+}(PG_{n})\big)}}\leq a\bigg)-{\int}^{a}_{-\infty}\frac{1}{\sqrt{2\pi}}e^{-\frac{t^{2}}{2}}dt \mid =0.
\end{eqnarray*}
\end{thm}

\noindent{\bf Proof of Theorem 3.4.}  Obtained by formula (2.20) yields, one sees that
\begin{eqnarray*}
Kf^{+}(PG_{n+1})=Kf^{+}(PG_{n})+\sum_{v\in V_{PG_{n}}}\big(12+5d(v)\big)r(u_{n},v)+203n+35.
\end{eqnarray*}
 Let $I_{n}=\sum\limits_{v \in V_{PG_{n}}}\big(12+5d(v)\big)r(u_{n},v)$. Then, we obtain
\begin{eqnarray}
Kf^{+}(PG_{n+1})=Kf^{+}(PG_{n})+I_{n}+203n+35.
\end{eqnarray}

Recalling that $Z_{n}^{1}$ and $Z_{n}^{2}$ are random variables, we obtain the following Facts.

\noindent{\bf Fact 3.4.1.} $I_{n}Z_{n}^{1}=(I_{n-1}+216n-133)Z_{n}^{1}.$

\noindent{\bf Proof.} If $Z_{n}^{1}$ = 0, the result is distinct. Then, we just take into account $Z_{n}^{1}$ = 1, which indicates $PG_{n} \rightarrow PG^{1}_{n+1}$. In this fact, $u_{n}$ (of $PG_{n}$) coincides with $x_{2}$ or $x_{5}$ (of $O_{n}$), see Fig. 3. In this situation, $I_{n}$ becomes
\begin{align*}
\begin{split}
&\sum_{v \in V_{PG_{n}}}\big(12+5d(v)\big)r(x_{2},v)\\
&=\sum_{v \in V_{PG_{n-1}}}\big(12+5d(v)\big)r(x_{2},v)+\sum_{v \in V_{O_{n}}}\big(12+5d(v)\big)r(x_{2},v)\\
&=\sum_{v \in V_{PG_{n-1}}}\big(12+5d(v)\big)\big(r(u_{n-1},v)+1+\frac{4}{5}\big)+12\sum_{v \in V_{O_{n}}}r(x_{2},v)+8\sum_{v \in V_{O_{n}}}d(v)r(x_{2},v)\\
&=\sum_{v \in V_{PG_{n-1}}}\big(12+5d(v)\big)r(u_{n-1},v)+\frac{9}{5}\sum_{v \in V_{PG_{n-1}}}\big(12+5d(v)\big)+12\times 4+5\times \frac{44}{5}\\
&=Q_{n-1}+\frac{9}{5}\big(12\times5(n-1)+5\times(12n-13)\big)+48+44\\
&=Q_{n-1}+216n-133.
\end{split}
\end{align*}
Thus, we obtain the desired Fact.

\noindent{\bf Fact 3.4.2.} $I_{n}Z_{n}^{2}=(Q_{n-1}+264n-181)Z_{n}^{2}.$

Similar to the proof of Fact 4.2.1, we only consider the fact $Z_{n}^{2}=1$, that is $PG_{n} \rightarrow PG^{2}_{n+1}$. In the same way, we omit the details.

Noting that $Z_{n}^{1}+Z_{n}^{2}=1$, by the above discussions, it holds that
\begin{eqnarray*}
I_{n}&=&I_{n-1}+(216Z_{n}^{1}+264Z_{n}^{2})n-(133Z_{n}^{1}+181Z_{n}^{2})\\
&=&Q_{n-1}+n\tilde{U}_{n}-\tilde{V}_{n},
\end{eqnarray*}
where for each $n$,

~~~~~~~~~$\tilde{U}_{n}=216Z_{n}^{1}+264Z_{n}^{2}$,~~~~$\tilde{V}_{n}=133Z_{n}^{1}+181Z_{n}^{2}$.
\\\text{Therefore, by Eq. (3.29),}\par

\begin{align}
Kf^{+}(PG_{n})=& Kf^{+}(PG_{1})+\sum_{l=1}^{n-1}I_{l}+\sum_{l=1}^{n-1}(203l+35)\nonumber\\
=& Kf^{+}(PG_{1})+\sum_{l=1}^{n-1}\big(\sum_{m=1}^{l-1}(I_{m+1}-I_{m})+I_{1}\big)+\sum_{l=1}^{n-1}(203l+35)\nonumber\\
=& Kf^{+}(PG_{1})+\sum_{l=1}^{n-1}\sum_{m=1}^{l-1}\big((m+1)\tilde{U}_{m+1}-\tilde{V}_{m+1}\big)+(n-1)I_{1}+\sum_{l=1}^{n-1}(203l+35)\nonumber\\
=& Kf^{+}(PG_{1})+\sum_{l=1}^{n-1}\sum_{m=1}^{l-1}\big((m+1)\tilde{U}_{m+1}-\tilde{V}_{m+1}\big)+O(n^2)\nonumber.
\end{align}

Suppose that

~~~~~~~~~$Var(\tilde{U}_{m})=\sigma^{2}_{4}$,~~~~$Var(\tilde{V}_{m})=\tilde{\sigma}^{2}_{4}$,~~~~$Cov(\tilde{U}_{m},\tilde{V}_{m})=r_{4}$.

The rest of proof is similar to the above Theorem, and the details are omitted here.\hfill\rule{1ex}{1ex}\

\section{Conclusion}
In this paper, we obtained the expected values and variances of Gutman index, Schultz index, multiplicative degree-Kirchhoff index and additive degree-Kirchhoff index about a class of the random chain networks. It was calculated and observed that under the same conditions, the expected value of Gutman index was the largest and that of multiplicative degree-Kirchhoff index was the lowest. Meanwhile, we found that they all approximately obeyed the normal distribution.


\begin{thebibliography}{99}
\small \setlength{\itemsep}{-.8mm}
\bibitem{1} M.E.J. Newman, $Networks$ (Oxford University, New York, 2018).

\bibitem{2} J.A. Bondy, U.S.R. Murty, $Graph\ Theory$, Springer, New York, 2008.

\bibitem{3} H. Wiener, Structrual determination of paraffin boiling points, $J.\ Am.\ Chem.\ Soc.$ 69 (1947):17-20.

\bibitem{4} R.C. Entringer, D.E. Jackson, and D.A. Snyder, Distance in graphs, $Czechoslovak\ Math.\ J.$ 26 (1976):283-296.

\bibitem{5} S. Mukwembi, and S. Munyira, MunyiraDegree distance and minimum degree, $Bull.\ Aust.\ Math.\ Soc.$ 87(2013):255-271.

\bibitem{6} H. Chen, and F. Zhang, Resistance distance and the normalized laplacian spectrum, $Discrete\ Appl.\ Math.$ 155, no. 5 (2007):654-661.

\bibitem{7} W.C. Shiu, P.C.B. Lam, The Wiener number of the hexagonal net, $Discrete\ Appl.\ Math.$ 73 (1997) 101–111.

\bibitem{8} W.C. Shiu, P.C.B. Lam, Wiener number of pericondensed benzenoidmolecule systems, $Congr.\ Numer.$ 126 (1997) 113–124.

\bibitem{9} W.C. Shiu, P.C.B. Lam, I. Gutman, Wiener number of hexagonal bitrapeziums and trapeziums, $Bull. Acad. Serbe Sci. Arts (Cl.\ Sci.\ Math.\ Nat.)$ 114 (1997) 9–25.

\bibitem{10} W.C. Shiu, C.S. Tong, P.C.B. Lam, Wiener number of hexagonal jagged-rectangles, $Discrete\ Appl.\ Math.$ 80 (1997) 83–96.

\bibitem{11} W.C. Shiu, P.C.B. Lam, K.K. Poon, On Wiener numbers of polygonal nets, $Discrete\ Appl.\ Math.$ 122 (2002) 251–261.

\bibitem{12} H. Deng, Wiener indices of spiro and polyphenyl hexagonal chains, $Math.\ Comput.\ Modelling$ 55 (2012) 634–644.

\bibitem{13} A.A. Dobrynin, R. Entringer, I. Gutman, Wiener index of trees: theory and applications, $Acta\ Appl.\ Math.$ 66 (2001) 211–249.

\bibitem{14} A.A. Dobrynin, I. Gutman, S. Klav\v{z}ar, P. \v{Z}igert, Wiener index of hexagonal systems, $Acta\ Appl.\ Math.$ 72 (2002) 247–294.

\bibitem{15} W. Yang, F. Zhang, Wiener index in random polyphenyl chains, $MATCH\ Commun.\ Math.\ Comput.\ Chem.$ 68 (2012) 371–376.

\bibitem{16} D.J. Klein, and M. Randi$\acute{c}$, Resistance distance, $J.\ Math.\ Chem.$ 12, no. 1 (1993):81-95.

\bibitem{17} A. Georgakopoulos, Uniqueness of electrical currents in a network of finite total resistance, $J.\ Lond.\ Math.\ Soc.$ 82, no. 1 (2010):256-272.

\bibitem{18} D. Bonchev, A. T. Balaban, X. Liu, and D. J. Klein, Molecular cyclicity and centricity of polycyclic graphs. I. Cyclicity based on resistance distances or reciprocal distances, $International\ Journal\ of\ Quantum\ Chemistry,$ vol. 50, no. 1, pp. 1–20, 1994.

\bibitem{19} D. Babi\'{c}, D. J. Klein, I. Lukovits, S. Nikoli\'{c}, and N. Trinajsti\'{c}, Resistance distance matrix: a computational algorithm and its application, $International\ Journal\ of\ Quantum\ Chemistry,$ vol. 90, pp. 166–176, 2002.

\bibitem{20} J. L. Palacios, Closed-form formulas for Kirchhoff index, $International\ Journal\ of\ Quantum\ Chemistry,$ vol. 81, no. 2, pp. 135–140, 2001.

\bibitem{21} E. Estrada and N. Hatano, Topological atomic displacements, Kirchhoff and Wiener indices of molecules, $Chemical\ Physics\ Letters,$ vol. 486, no. 4–6, pp. 166–170, 2010.

\bibitem{22} J. Huang, S. Li, and L. Sun, The normalized Laplacians, degree-Kirchhoff index and the spanning trees of linear hexagonal chains, $Discrete\ Applied\ Mathematics,$ vol. 207, pp. 67–79, 2016.

\bibitem{23} Y. Yang and H. Zhang, Kirchhoff index of linear hexagonal chains, $International\ Journal\ of\ Quantum\ Chemistry,$ vol. 108, no. 3, pp. 503–512, 2008.

\bibitem{24} Q. Deng and H. Chen, On extremal bipartite unicyclic graphs, $Linear\ Algebra\ and\ Its\ Applications,$ vol. 444, pp. 89–99, 2014.

\bibitem{25} L. Feng, G. Yu, K. Xu, and Z. Jiang, A note on the Kirchhoff index of bicyclic graphs, $Ars\ Combinatoria,$ vol. 114, pp. 33–40, 2014.

\bibitem{26} X. Gao, Y. Luo, and W. Liu, Resistance distances and the Kirchhoff index in Cayley graphs, $Discrete\ Applied\ Mathematics,$ vol. 159, no. 17, pp. 2050–2057, 2011.

\bibitem{27} J. Huang, S. Li, and Q. Zhao, On extremal bipartite bicyclic graphs, $Journal\ of\ Mathematical\ Analysis\ and\ Applications,$ vol. 436, no. 2, pp. 1242–1255, 2016.

\bibitem{28} J. Huang, S. Li, and X. Li, The normalized Laplacian, degree-Kirchhoff index and spanning trees of the linear polyomino chains, $Applied\ Mathematics\ and\ Computation,$ vol. 289, pp. 324–334, 2016.

\bibitem{29} J.-B. Liu, X.-F. Pan, L. Yu, and D. Li, Complete characterization of bicyclic graphs with minimal Kirchhoff index, $Discrete\ Applied\ Mathematics,$ vol. 200, pp. 95–107, 2016.

\bibitem{30} J. L. Palacios, On the Kirchhoff index of graphs with diameter 2, $Discrete\ Applied\ Mathematics,$ vol. 184, pp. 196–201, 2015.

\bibitem{31} M.R. Farahani, Hosoya, Schultz, Modified Schultz polynomials and their topological indices of benzene molecules: First members of polycyclic aromatic hydrocarbons (PAHs), $Int.\ J.\ Theor.\ Chem.$ 1, no. 2 (2013):9-16.

\bibitem{32} I. Gutman, Selected properties of the Schultz molecular topological index, $J.\ Chem.\ Inf.\ Comput.\ Sci.$ 34 (1994):1087-1089.

\bibitem{33} A. Heydari, On the modified Schultz index of $C_{4}C_{8}(S)$ nanotubes and nanotorus, $Digest.\ J.\ Nanomater.\ Biostruct.$ 5, no. 1 (2010):51-56.

\bibitem{34} S. Mukwembi, and S. Munyira, Munyira degree distance and minimum degree, $Bull.\ Aust.\ Math.\ Soc.$ 87 (2013):255-271.

\bibitem{35} H. Chen, and F. Zhang, Resistance distance and the normalized laplacian spectrum, $Discrete\ Appl.\ Math.$ 155, no. 5 (2007):654-661.

\bibitem{36} I. Gutman, L. Feng, and G. Yu, Degree resistance distance of unicyclic graphs, $Trans.\ Comb.$ 1, no. 2 (2012):27-40.

\bibitem{37} S.B. Huang, J. Zhou, and C.J. Bu, Some results on Kirchhoff index and degree-Kirchhoff index, $MATCH\ Commun.\ Math.\ Comput.\ Chem.$ 75, no. 1 (2016):207-222.

\bibitem{38} M. Somodi, On the Ihara zeta function and resistance distance-based indices, $Linear\ Algebra\ Appl.$ 513 (2017):201-209.

\bibitem{39} J.L. Zhang, X.H. Peng, and H.L. Chen, The limiting behaviours for the Gutman index, Schultz Index, multiplicative degree-Kirchhoff index and additive degree-Kirchhoff index of a random polyphenylene chain, $Discrete\ Appl.\ Math.$ 299 (2021):62-73.

\bibitem{40} W. Yang, and F. Zhang, Wiener index in random polyphenyl chains, $MATCH\ Commun.\ Math.$ $Comput.\ Chem.$ 68 (2012):371-376.

\bibitem{41} L. Ma, H. Bian, B.J. Liu, and H.Z. Yu, The expected values of the Wiener indices in the random phenylene and spiro chains, $Ars\ Combinatoria$ 130 (2017): 267-274.

\bibitem{42} G.H. Huang, M.J. Kuang, and H.Y. Deng, The expected values of Kirchhoff indices in the random polyphenyl and spiro chains, $Ars\ Mathematics\ Contematica$ 9, no. 1 (2014):197-207.

\bibitem{43} S. Wei, and W.C. Shiu, Enumeration of Wiener indices in random polygonal chains, $J.\ Math.\ Anal.\ Appl.$ 469, no. 2 (2019):537-548.

\bibitem{44} L.L. Zhang, Q.S. Li, S.C. Li, and M.J. Zhang, The expected values for the Schultz index, Gutman index, multiplicative degree-Kirchhoff index and additive degree-Kirchhoff index of a random polyphenylene chain, $Discrete\ Applied\ Mathematics$ 282 (2020):243-256.

\bibitem{45} J.L. Zhang, X.H. Peng, and H.L. Chen, The limiting behaviours for the Gutman index, Schultz index, multiplicative degree-Kirchhoff index and additive degree-Kirchhoff index of a random polyphenylene chain, $Discrete\ Appl.\ Math.$ 299 (2021):62-73.

\bibitem{46} H.C. Liu, R.W. Wu, L.H. You, Three types of Kirchhoff indices in the random cyclooctane Chains, $Journal\ of\ South\ China\ Normal\ University\ ( Natural\ Science\ Edition)$ 53, no. 2 (2021):96-103.

\bibitem{47} S. Wei, X. Ke, Y. Wang, Wiener indices in random cyclooctane chains, $Wuhan\ Univ.\ J.\ Nat.\ Sci.$ (2018), submitted for publication.

\bibitem{48} G. Huang, M. Kuang, and H. Deng, The expected values of Kirchhoff indices in the random polyphenyl and spiro chains, $Ars\ Mathematica\ Contemporanea,$ vol. 9, pp. 197–207, 2015.

\bibitem{49} V.V. Petrov, Limit theorems of probability theory sequences of independent random variables, $Oxford\ University\ Press$, 1995.

\bibitem{50} W. Feller, An introduction to probability theory and its applications, $Wiley\ Series\ in\ Probability\ and\ Statistics$, 1971.




\end{thebibliography}
\end{document}